\documentclass[final,3p,times,authoryear]{elsarticle}

\usepackage{amssymb}

\PassOptionsToPackage{hyphens}{url}

\usepackage[dvipsnames,table]{xcolor}
\usepackage{fancyhdr,graphicx,amsmath,amssymb}
\usepackage{subfigure}
\usepackage{multirow}
\usepackage{rotating}
\usepackage{mathtools}

\usepackage[ruled,vlined,linesnumbered]{algorithm2e}
\usepackage{hyperref}
\usepackage{natbib}
\bibpunct[, ]{(}{)}{,}{a}{}{,}

\usepackage{tikz}
\usetikzlibrary{arrows.meta, positioning, calc, fit}
\usepackage{longtable}

\usepackage{mathtools}
\usepackage[normalem]{ulem}
\usepackage{pifont}
\graphicspath{{Figures/}} 
\hypersetup{colorlinks,%
citecolor={blue}, 
urlcolor={blue},
linkcolor={blue},
breaklinks={true}
}

\usepackage{makecell, multirow, rotating}
  
\newcommand{\vrp}{VRP-SDO}
\newcommand{\vrpx}{VRP-SD}
\newcommand{\vrpvcsd}{VRP-VCSD}
\newcommand{\vrppo}{VRP-PFCC}
\newcommand{\dqnco}{DQNCO}
\newcommand{\qnco}{SQNCO}
\newcommand{\irp}{IRP}
\newcommand{\igp}{IGP}
\newcommand{\ihp}{IHP}
\newcommand{\iqnco}{ISQNCO}
\newcommand{\idqnco}{IDQNCO}
\newcommand{\iqncot}{ISQNCO$^+$}
\newcommand{\idqncot}{IDQNCO$^+$}

\usepackage{marginnote}

\begin{document}

\begin{frontmatter}



\title{A Deep Reinforcement Learning Algorithm for the Vehicle Routing Problem with Stochastic Demands and Outsourcing}

\author[ets,gerad,cirrelt]{Mohsen Dastpak\texorpdfstring{\corref{cor1}}{}}
\ead{mohsen.dastpak.1@ens.etsmtl.ca}
\cortext[cor1]{Corresponding author}
\author[ets,gerad,cirrelt]{Fausto Errico}
\author[cirrelt,polimi]{Ola Jabali}
\affiliation[ets]{organization={Department de génie de la construction, École de technologie supérieure},
            addressline={}, 
            city={Montréal},
            postcode={H3C 1K3}, 
            state={QC},
            country={Canada}}
\affiliation[gerad]{organization={GERAD},
            addressline={}, 
            city={Montréal},
            postcode={H3T 1J4}, 
            state={QC},
            country={Canada}}
\affiliation[cirrelt]{organization={CIRRELT},
            addressline={}, 
            city={Montréal},
            postcode={H3C 3J7}, 
            state={QC},
            country={Canada}}
\affiliation[polimi]{organization={Dipartimento di Elettronica, Informazione e Bioingegneria, Politecnico di Milano},
            addressline={Piazza Leonardo da Vinci 32}, 
            city={Milano},
            postcode={20133}, 
            state={},
            country={Italy}}

\begin{abstract}
We introduce the vehicle routing problem with stochastic demands and outsourcing options (\vrp{}). In this problem,  we consider a logistics service provider  (LSP) that must partition a set of customer requests with uncertain demands into customers outsourced to a common carrier and customers committed to its own fleet. The former are exclusively served by the common carrier, whereas the LSP fully serves its committed customers with a fixed fleet. The latter induces a vehicle routing problem with stochastic demand (\vrpx{}), which is solved dynamically while allowing multiple vehicle visits to a customer. Specifically, the demand of  committed customers is revealed once visited. If the vehicle capacity is not sufficient to serve it, the residual demand may be served by other vehicles or by the same vehicle after having restocked at the depot. Moreover, vehicles are allowed to replenishment trips to the depot. Vehicles exceeding regular work-shift duration incur an overtime cost. The unit outsourcing cost is inversely correlated with the expected total demand of the outsourced customers. The \vrp{} objective is to minimize expected total travel and overtime costs, as well as  the outsourcing costs. 

We propose an iterative two-level solution methodology for the \vrp{}. The first level partitions customers into committed and outsourced subsets, while the second level  solves the resulting \vrpx{} and estimates its expected routing cost. Solving this stochastic dynamic problem for a fixed set of committed customers is computationally demanding. Given that the \vrpx{}  must be evaluated repeatedly for different subsets of committed customers throughout the iterative search, solving it from scratch at every iteration is computationally impractical. Therefore, we propose learning  an offline routing policy which estimates routing costs almost instantaneously for any subset of committed customers. We demonstrate our methodology using an iterated local search (ILS) algorithm that primarily establishes the first-level partitions.
    
We formulate the second level problem as a Markov decision process and solve it using a deep Q-network (DQN). To represent the state of the system, we propose a graph attention network (GAT) that aggregates customer and vehicle information based on their relevance to the vehicle selecting an action, capturing the spatial structure induced by their relative locations. The DQN is trained offline on instances with variable customer cardinality and location, yielding a routing policy applicable to any daily realization of customers. We further introduce an online fine-tuning step that improves the accuracy of the routing cost approximation during the search.
    
Extensive experiments show that  our DQN for the \vrpx{} reduces routing costs by 19.6\% relative to a state-of-the-art methodology, and by at least 29.6\% when compared to classical routing heuristics. Our overall solution algorithm for the \vrp{} yields savings of 13.7\% on average relative to the version whose routing policy does not use  the GAT representation. Furthermore, our method generates the high-quality decisions within minutes, whereas the benchmarks that do not use an offline-trained cost estimator require more than one hour on average.

\end{abstract}


\end{frontmatter}

\section{Introduction}\label{sec:intro}
The volume of last mile delivery services has significantly increased over the past years. In 2024 the number of packages delivered in the US was 22.4 billion, up from 15.4 billion in 2019 \citep{PitneyBowes2025}. 
Logistics service providers (LSPs) operating next-day delivery services typically manage a fixed set of resources (e.g., trucks).
These are often insufficient to handle varying demand volumes. Furthermore, expanding resources to serve areas away from an LSP's center of activity is often economically unviable.  As a result,  many LSPs structurally embed outsourcing services to fully satisfy the demand of their customers. 
For example, while Amazon has significantly grown its internal logistics network, analyses indicate that third-party carriers continue to handle a substantial number  of packages, estimated at nearly 30\% of total volume as of late 2021 \citep{Tarasov2021}.

The use of outsourcing services by LSPs is governed by  contractual agreements. In the context of next-day delivery, operational practices described in the literature indicate that an LSP is often obliged to communicate its \textit{outsourced customers} of a given day by the morning of that day (see \cite{CanadaPostParcelGuide} for the \emph{Parcel Services Customer Guide} by Canada Post). In turn, the LSP commits to serving the remainder set of customers, which we refer to as \textit{committed customers}. This setting has been treated in the scientific literature as the VRP with private fleet and common carrier (\vrppo{}) \citep[e.g.,]{Gahm2017, Baller2020}. As most routing problems, the \vrppo{} is subject to uncertainty on multiple dimensions in practice (e.g., travel times and customer demands).

Routing under uncertainty may be cast under multiple modeling frameworks, such as stochastic programming, robust optimization, and stochastic dynamic programming (SDP), see  \cite{Oyola2017, Oyola2018}. The latter has been gaining increased attention, as LSPs are gradually accepting to take decisions dynamically. Moreover, the methodological advances in stochastic dynamic routing problems have made the theoretical advantages of such models a tangible reality. Such advances are largely attributed to the integration of reinforcement learning in SDP problems. In this paper, we consider a stochastic dynamic variant of the \vrppo{}, where customer demands are uncertain. Accordingly, we introduce the vehicle routing problem with stochastic demands and outsourcing options (\vrp{}). In this problem, the LSP receives customer requests at the beginning of the operational day, requiring a rapid decision regarding which customers to outsource to a third-party carrier and which to commit to. As in \cite{Dabia2019}, we consider that  the cost per unit of outsourcing is inversely correlated with the expected total demand of the outsourced customers. 
Furthermore, the demand of customers is uncertain, and we assume that outsourced customers are exclusively served by the common carrier, whereas the LSP fully serves its committed customers, while possibly incurring overtime costs. In particular, we consider that vehicles are dynamically assigned to customers, and allow different  vehicles to visit a customer, if a single visit is not sufficient.

The LSP must solve the \vrp{} every day within a limited computing budget. From an application standpoint, the major limitation to the computing budget is the need to communicate the set of outsourced customers early in the day. The cost of partitioning a set of customers between outsourced ones and committed ones is not trivial to compute. In particular,  the set of committed customers induces a stochastic dynamic problem, which we denote as the vehicle routing problem with stochastic demands (VRP-SD). 
Therefore, we propose an iterative two-level solution methodology. Such a hierarchy is common in the deterministic \vrppo{} literature \citep{Stenger2013, Vidal2016, Dabia2019}. Specifically, the first-level partitions the customers into the two groups, allowing the immediate cost evaluation of the outsourced customers. In order to maintain a viable iterative procedure between the two stages, evaluating the costs of various sets of committed customers needs to be efficiently done. Indeed, any algorithm based on the exploration of the first-level decision space would require many evaluations of the second-level decision problem. However, solving the second-level routing problem (\vrpx{}) for any fixed set of committed customers is computationally intensive. This subproblem, especially due to its stochastic and dynamic nature, may require seconds to several minutes  for its solution \citep{Erera2010, Mendoza2016}. While these run times may seem acceptable in isolation, they become a significant bottleneck when the problem must be solved repeatedly during the execution of the solution algorithm. The strategy proposed in this paper aims then at transferring a consistent part of the computational burden from everyday online execution to an offline learning approach.

The choice of an offline learning approach is motivated by the fact that the complete set of customers to serve in the \vrp{} differs from day to day. Furthermore, this complete set of customers is typically known shortly prior to establishing the sets of committed and outsourced customers. Thus, it is impractical to retrain a reinforcement learning algorithm on a new set of customers every day.
As our methodology relies on deploying a trained algorithm, we propose training the algorithm on a  VRP-SD while also considering the set of customers to be variable in terms of cardinality and  location. This allows estimating the routing cost of a \vrpx{} in instant time for any daily realization of customers. We showcase our methodology by using an iterated local search (ILS) algorithm that  primarily  establishes the first-level partitions. However, it is worth mentioning that our methodology can be applied to other metaheuristics. 
 
We formulate the \vrpx{} as a Markov decision process (MDP), whose value function approximates the routing costs associated with a given set of committed customers as an input. 
To solve the resulting MDP, we build upon the methodology put forward by \cite{Dastpak2023}. The uncertain parameters considered in that paper resemble those of the \vrpx{}. However, \cite{Dastpak2023} deal with maximizing the expected served demands while meeting vehicle capacity and time restrictions. One of the main methodological contributions of that paper lies in using an observation function to guide decision-making.
The proposed observation function represents the state of the system as a fixed-size vector, comprising a limited yet  relevant  set of features. These features were hand-selected based on various experiments. Moreover, the observation function heuristically chooses a set of target customers, which may be visited. As such, the proposed observation function does not  capture the underlying graph structure induced by customer locations and therefore fails to represent spatial relationships such as proximity and neighborhood interactions. 
Lastly, customer information is discretized using a grid-based technique, which may greatly impact the precision of the overall methodology. To address these limitations, we propose a graph attention network (GAT) architecture to represent the  observation function. 
This architecture employs an attention mechanism \citep{Vaswani2017} to aggregate customers' and vehicles' information into a fixed-size vector by considering their relevance to the vehicle that needs to take an action. One of the core contributions of this paper is to demonstrate the power of GATs to represent a concise but informative observation vector.

Considering the GAT, we propose a deep Q-network (DQN) for the \vrpx{}. DQNs are reinforcement learning  algorithms combining Q-learning with deep neural networks. We propose a Q-learning algorithm using a two-layer neural network to approximate the value of observation-action pairs, also referred to as the Q-network. In particular, we train the Q-network alongside the proposed GAT architecture of the observation function, resulting in a deep Q-network. We denote the proposed routing method for solving the \vrpx{} as \dqnco{}. 
We employ several techniques to enhance the overall performance of the Q-learning algorithm,
including replay memory~\citep{Mnih2015} and double Q-network~\citep{VanHasselt2016}. 
We name the entire solution method, which is the ILS algorithm with the embedded routing policy, as \idqnco{}. To further enhance the performance of the \idqnco{}, we introduce an online fine-tuning step inside the ILS algorithm to improve the accuracy of the routing cost value function. The \idqnco{} with this additional step is denoted as \idqncot{}.

We conduct comprehensive experiments to evaluate the performance of the proposed algorithm, comparing it with six benchmark methods. 
We first show the superiority of the GAT-based architecture introduced in the \dqnco{} over all the other benchmarks, including \qnco{}, the routing policy from \cite{Dastpak2023}, based on the hand-engineered observation function. 
We then demonstrate that \idqncot{} significantly outperforms all the other benchmarks. Notably, our method reduces overall costs by more than 22.3\% compared to benchmarks with basic routing policies and by an average of 13.7\% compared to \iqncot{}, which combines the same ILS algorithm and fine-tuning step with the \qnco{} routing policy.  
Summarizing, the contributions of this paper are as follows:
\begin{itemize}
    \item We introduce the \vrp{}, as a stochastic  variant of the \vrppo{}, and  propose a two-level formulation for it.
    \item We formulate the outsourcing problem as a binary integer programming problem and the \vrpx{} as an MDP. 
    \item We propose a DQN-based algorithm that is trained offline and rapidly estimates the routing costs of any committed customer set within the ILS.
    \item We propose a GAT architecture to represent the active vehicle's observation, taking into account the graph-like structure of customers and vehicles in the service area, as well as their relevance to the active vehicle.
    \item We conduct extensive experiments and demonstrate that our method outperforms six benchmarks. 
\end{itemize}

The rest of the paper is organized as follows. We survey the literature around our problem in Section~\ref{sec:literature}.
In Section~\ref{sec:problem}, we describe the \vrp{} and provide its formulation. In Section~\ref{sec:method} we describe the proposed solution methodology, and we present computational results in Section~\ref{sec:experiments}. Finally, we draw our conclusions in Section~\ref{sec:conclusion}.

\section{Literature Review}\label{sec:literature}
In this section we provide an overview of the scientific literature related to the  \vrp{}  and to our chosen methodology. In Section \ref{sec:relatedproblems}, we position  the \vrp{} within the landscape of relevant problems. In Section \ref{sec:relatedmethods} we survey the main solution methods pertaining to closely related variants of the \vrp{} and position our methodological contribution.

\subsection{Related problems}\label{sec:relatedproblems}
The \vrp{} encompasses two main complicating dimensions: stochastic customer demands and the determination of outsourced customers. With the latter being an option prior to commencing the routes and revealing any customer demand. As the combination of both dimensions is relatively rare in the literature, we survey the relevant contributions of each in separation.

Uncertain demands have been widely explored in the routing literature  \citep[see][for a comprehensive survey]{Parada2024}, where it is typically assumed that  customer demand is revealed upon the first  arrival of a vehicle. As a result,  if the available vehicle capacity is not sufficient  to satisfy customer demand,  a  \textit{classical} recourse action is predominantly applied, whereby the vehicle restocks its capacity at the depot and returns to the customer \citep{Laporte2002, Jabali2014}. This recourse action  is often complemented by allowing vehicles to perform preventive returns to the depot (e.g., \cite{Louveaux2018, Florio2020}).

The problem settings using classical recourse, as well as preventive returns,  implicitly assume that each customer is served exclusively by a single vehicle. Thus, any recourse action is performed by the vehicle that had visited the customer for the first time. This assumption facilitates modeling the resulting problem as a two-stage stochastic program, where the first stage establishes fixed routes for the vehicles, which are maintained irrespectively of the outcome of the second stage, thus having second stage optimization problems that are separable by vehicle (e.g., \cite{Hoogendoorn2023}). 
Such  unique vehicle-to-customer assignment assumptions  have been relaxed in the literature (e.g., \cite{Goodson2016}) by assuming that all  vehicles may  serve  a customer in additional visits, if needed. This flexibility often comes in the form of dynamically establishing the customer visits for each vehicle.  Considering the committed customers, we  dynamically allocate vehicles to customers while allowing customers to be visited by multiple vehicles, in case a single  visit is not sufficient.  We consider classical recourse as well as preventive returns.

Irrespectively of the adopted recourse strategy, the variability in customer demand leads to  variability in route duration. This introduces challenges  when vehicles are constrained by regular working hours. For example, \cite{Erera2010, Goodson2016}, and \cite{Dastpak2023} assume hard duration limits for each vehicle trip. \cite{Erera2010} consider a variable number of vehicles to ensure serving all customers under the duration limit constraints, while the latter two works maintain a fixed set of vehicles and aim to maximize the total served demand within the duration limit. In \vrp{}, similar to \cite{Mendoza2016}, we consider a fixed set of available vehicles with a soft duration limit, where exceeding regular working hours is possible at an additional cost.

A common approach for controlling costs in routing problems with limited resources, such as fleet size, vehicle capacity, and trip duration, involves outsourcing the service of a subset of customers. This is often addressed as a deterministic problem, which is referred to as the VRP with private fleet and common carrier (\vrppo{}). For instance, \cite{Ceschia2011, Moon2012, Brito2015, Gahm2017} investigate cases in which the LSP can rent vehicles in addition to its own fleet. In these cases, a daily cost or a distance-dependent cost for rented vehicles is optimized.  Alternatively, \cite{Chu2005, Bolduc2008,Cote2009}, and \cite{Poon2022} study problems in which  outsourced customers have a given outsourcing cost.
In contrast, \cite{Baller2020} and \cite{Zhu2022} consider a cost per unit of outsourced demand. In particular, \cite{Baller2020} consider outsourcing a portion of a customer's demand.
In real-world applications, common carrier companies may offer volume-based discounts based on either individual customer demands \citep{Stenger2013} or aggregated volume of outsourced demands \citep{Dabia2019}. Similar to the latter, we consider that outsourcing rates decrease as the total quantity of outsourced demands increases. Given that in the \vrp{} outsourcing decisions are to be established prior to revealing the customer demand values, the considered outsourcing costs are based on the total expected demand values. Furthermore, we acknowledge that the \vrp{} shares modeling features with the capacitated profitable routing problem \citep{Archetti2009}, the team orienteering problem \citep{Boussier2007, Vidal2016}, and the prize-collecting problem \citep{Moradi2025}.

The previously mentioned \vrppo{}s deal with deterministic problems. To the best of our knowledge, only \cite{Zhu2022} address a dynamic VRP with stochastic demands where outsourcing is performed when a vehicle cannot fully serve a customer. Therefore, outsourcing is not a primary assignment  decision, but is rather a recourse penalty per unit.
However, in the \vrp{}, the outsourced customers are determined at the beginning of the day, which is aligned with common practice. To the best of our knowledge, no prior work jointly addresses vehicle routing  with demand uncertainty where outsourcing decisions are made a priori.  Indeed, the \vrp{} combines complex  modeling features from two strands of the literature. Vehicles are dynamically dispatched to customers, thus impeding the separation of subproblem on a per vehicle basis. Moreover, the outsourcing costs are a function of the total expected demand value of outsourced customers, which creates significant dependencies among these decisions.

\subsection{Solution methods}\label{sec:relatedmethods}
The solution approaches to the \vrppo{} vary depending on the type of outsourcing service considered by the LSP.
In problems involving renting vehicles, a common approach is to establish routes for both the private and rental fleets \citep{Brito2015, Sawadsitang2018, Gahm2017}. Consequently, the outsourcing decisions become an integral part of the routing process. 
In problems where the LSP chooses a subset of customers to outsource and does not route them, a virtual vehicle route is typically used to model  the common carrier \citep{Stenger2013}.
 Specifically, in these cases,  solution methods
often  assign and route a subset of customers to $m$ private vehicles, whereas the remaining subset is assigned to a virtual route representing the outsourcing decisions.

\cite{Goeke2019, Dabia2019, Su2023}, and \cite{Moradi2025} have proposed branch-price-and-cut algorithms to address deterministic \vrppo{} variants. 
Additionally, several authors have proposed heuristic methods for the \vrppo{}. For example, \cite{Chu2005} consider outsourcing with a given cost per customer, employing a savings-based constructive procedure and local improvement. \cite{Bolduc2008} introduce a heuristic for the same problem, integrating a local descent on various neighborhood structures with a perturbation mechanism involving the swapping of customer pairs, and an improvement procedure. \cite{Cote2009} propose a tabu search to address the same problem and achieved improvements over previous results. Additionally, \cite{Potvin2011} present an enhanced tabu search heuristic with a neighborhood structure based on ejection chains. To address a multi-depot extension of the \vrppo{}, \cite{Stenger2013} proposed an adaptive variable neighborhood search, incorporating an adaptive mechanism for routing and customer selection within the shaking phase. Subsequently, \cite{Vidal2016} introduce  large neighborhoods for this problem, resulting in higher-quality solutions than  those of \cite{Stenger2013}. Considering a \vrppo{} with multiple common carrier types, \cite{Gahm2017} explored variants of variable neighborhood search.

The previously mentioned methodologies address deterministic problems. In the \vrp{} the customer demands are stochastic and committed customers must be fully served by the LSP fleet. Our solution methodology relies on an ILS algorithm, which is commonly used for routing problems, e.g., \citep{Maximo2024}. At each iteration, our methodology first determines the committed and outsourced customers, where evaluating the cost of the former is a vehicle routing problem with
stochastic demands. We model the \vrpx{}  as an MDP considering classical recourse as well as preventive returns, while allowing  multiple vehicles to serve the same customer. 
Two primary methodologies have been adopted to solve  MDP-based SVRP. Approximate policy iteration methods (e.g.,  policy gradient), which aim at directly developing the desired policy by iteratively improving an initial policy \citep{Li2021}. 
Conversely, approximate value iteration methods, such as Q-learning, compute a value function for state-action values (e.g., \cite{Li2021a, Chen2022}). We adopt a Q-learning algorithm, where we develop a  value function estimating routing costs.

MDP formulations are typically addressed using dynamic programming algorithms, yet these algorithms are susceptible to the curse of dimensionality.  This challenge is prevalent in stochastic routing problems, which are characterized by expansive combinatorial state and action spaces. Various mechanisms have been proposed to address these complexities, such as decomposing the problem into multiple single-vehicle subproblems \citep{Fan2006, Goodson2013}, or decentralizing (i.e., computing value functions for individual vehicles rather than the joint fleet) while either  neglecting \citep{Oda2018, OroojlooyJadid2019} or considering \citep{Li2019, Chen2019, Kullman2020} inter-vehicle collaboration.  
For a comprehensive review, we direct readers to \cite{Dastpak2023}. In this paper, to alleviate the  issues associated with large state and action spaces, we adopt two methodological  features  from their work. Namely, we use an  MDP with consecutive action selection and an observation function (MDP-CO). Consecutive action selection  entails that  only one vehicle makes a decision at each decision epoch, whereas the entire state of the system is aggregated into the  vehicle's observation. Specifically, \cite{Dastpak2023} proposed a simple observation function based on grid-based aggregation of all customers \citep[see, for example][]{Chen2019,Kullman2020} and a basis function representation \citep[see, for example][]{Joe2020} of \textit{target customers}, who are heuristically identified as promising candidates for the active vehicle. 

In this paper, we improve the observation function introduced by \cite{Dastpak2023} in two major and complementary ways, both realized within a single graph attention network (GAT). 
The first improvement concerns how node information is represented. In our setting, the nodes of the graph are the customers, each described by its own raw features (e.g., location and expected demand). \cite{Dastpak2023} treat these nodes in isolation and thereby overlook the spatial dependencies and topological relationships among them. Such relationships are naturally captured by graph neural networks (GNNs), which represent a node by aggregating information from its neighbors \citep{Scarselli2009} and have proven effective across a range of combinatorial optimization problems \citep{Cappart2021}.
An early example of such a GNN is the ``structure2vec'', which represents each node as a learned combination of its neighborhood features (e.g., coordinates and demand) and has been applied to the traveling salesman problem (TSP) \citep{Dai2016} and the VRP \citep{Lin2020}.
A more advanced way to perform such neighborhood aggregation is through attention mechanisms, which learn  the relevance of each neighbor to the node being represented.
This mechanism was first introduced by \cite{Bahdanau2015} to help long short-term memory (i.e., LSTM) and recurrent neural network (i.e., RNN)-based encoder-decoders to cope with the long-range sequences in natural language processing, as later scaled up in systems such as Google's neural machine translation \citep{Wu2016}.
Although originally designed for sequences, the underlying principle of attention (i.e., dynamically weighting the importance of input elements) transfers naturally to graph-structured data. 
The transformer encoder \citep{Vaswani2017} is the architecture that applies this principle to a set of nodes, encoding each one as a weighted combination of all the others.
In a vehicle routing context, transformer encoders have been employed to generate representations that capture each customer's relative position within the graph \citep{Xu2022}.
We build our GAT from such transformer-style attention layers. Accordingly, each node is encoded as a weighted combination of all other nodes, where the weights express the learned relevance (or alignment) between node pairs.

The proposed GAT addresses the second limitation of \cite{Dastpak2023}, which concerns how the node representations are reduced to the fixed-size observation of the active vehicle.
To represent a varying number of customers into a fixed-size vector, \cite{Dastpak2023} used a heuristic rule to select a subset of promising customers of fixed cardinality (i.e., target customers) in the observation. By retaining  only a subset of customers, this rule produces an incomplete  observation that may omit  important  information for  routing decisions.
We instead let the GAT decide, through another attention mechanism, how to represent the information of all customers and vehicles into a single fixed-size observation. Specifically, the GAT assigns every customer and vehicle a weight that reflects its relevance to the active vehicle, and then combines all of them into a single fixed-size vector. 
This way of summarizing a variable set of nodes into one fixed-size vector through attention follows \cite{Kool2019}, whose encoder-decoder architecture consolidates the customer set into a single vector before constructing a TSP tour, and \cite{Nazari2018}, who extend the idea to the VRP.
The approaches of \cite{Bono2020} and \cite{Zhang2020} are the most closely related to ours, as they also address multi-vehicle VRPs and embed customers and vehicles jointly through attention mechanisms. 
However, in both works the attention layers directly parameterize a routing policy by scoring customers for the next visit, whereas our GAT builds a fixed-size observation on which a value function is learned to estimate routing costs.

In conclusion, to the best of our knowledge, no existing solution method in the literature can be trivially adapted to the \vrp{}. We present a novel solution method that  decomposes the problem into outsourcing and routing components. The latter component is modeled as an MDP and solved  with a deep Q-network based algorithm. This  yields a fast and reliable mechanism for estimating routing costs, thus allowing an iterative evaluation of a large number of decomposed solutions. The pillar of this mechanism is a novel GAT structure for representing the system's state. Furthermore, we propose an online fine-tuning step inside the local search heuristic to further enhance the performance of the developed value function in real-time.

\section{The Vehicle Routing Problem with Stochastic Demands and Outsourcing Options} \label{sec:problem}
In Section \ref{sec:description}, we formally define the \vrp{}. In Section \ref{sec:formulation} we present a high-level mathematical formulation for it based on a two-level decision structure, namely, outsourcing decisions at the first level, and routing at the second (the \vrpx). In Section \ref{sec:MDP} we provide an MDP formulation of the \vrpx{}.

\subsection{Problem Description}\label{sec:description}
The LSP operates $m$ vehicles, denoted by $\mathcal{V}=\{1, 2, ..., m\}$, all initially located at the depot. We denote the depot and its location by $0$ and $l_0$, respectively.
We assume that vehicles are identical and have capacity $Q$.
Every morning, the LSP receives service requests from a set of customers $\mathcal{C} = \{1,\ldots,n\}$. The requests are to be served during the day. Thus, all service requests are known before routing and outsourcing decisions are made. Although pick-up and delivery interpretations are mathematically equivalent in our setting, we adopt the pick-up interpretation for consistency of exposition.
We assume that the set of customer requests  varies from one day to another.
Each customer $c\in\mathcal{C}$ is characterized by its location $l_c$ and a stochastic demand $d_c$, whose expected value $\bar{d}_c$ is assumed to be known. 
In particular, the demand of each customer $c\in\mathcal{C}$ is assumed to follow a known
probability distribution denoted $\Gamma^D_c = \Gamma^D(\bar{d}_c)$, where $\Gamma^D(\bar{d}_c)$ is a function that returns a probability distribution for a given input value $\bar{d}_c$.
We define a complete graph $G=(N, E)$ with the set of nodes $N=\{0\} \cup \mathcal{C}$ and the set of arcs $E=\{(i,j)~|~i,j \in N, ~ i\neq j\}$, which represents potential visit sequences between nodes.
The travel time between two nodes $i$ and $j$, denoted by $\tau_{ij}$, is defined for each $(i, j) \in E$.

At the beginning of the operating day, the LSP determines the subset of outsourced customers $\mathcal{C}^o\subset\mathcal{C}$. We let $\overline{\mathcal{C}}^o = \mathcal{C}\setminus \mathcal{C}^o$ denote the set of committed customers, which will be served by the LSP's fleet.
Partially outsourcing a customer's demand is not allowed, i.e., a customer's demand is either fully outsourced or fully served by the LSP.  
Similar to \cite{Dabia2019}, we assume an outsourcing cost $\Psi(\mathcal{C}^o)$ to be a piecewise linear function of the expected total outsourced demand $\sum_{c\in\mathcal{C}^o} \bar{d}_c$.  

The LSP's vehicles begin and end their routes at the depot. 
As customary in the literature, we assume that a customer's actual demand is observed when a vehicle visits the customer \citep{Florio2023}.
We allow the demand of each committed customer $c\in\overline{\mathcal{C}}^o$ to be collected over multiple visits. Specifically, if a vehicle is unable to completely serve a customer during a visit, it serves it to the fullest possible extent, and any remaining demand is collected by the same or other vehicles at a later time.
Vehicles are also allowed to perform a \textit{preventive unloading} action, meaning that they can visit the depot to unload even if their capacity has not been fully used yet.
Vehicle operation time is charged at one unit of cost per unit of time up to a regular work-shift duration $L$, and at an \textit{overtime} cost of $\phi$ per time unit beyond it.
We assume that a vehicle may remain at the depot and thus terminates its operation, provided that at least one vehicle remains operative to attend remaining customers.
The objective is to minimize the total expected cost, including outsourcing, travel, and overtime cost.

\subsection{A high-level formulation of the \vrp{}}\label{sec:formulation}
We now present a two-level mathematical formulation for the  \vrp{}. The first-level problem, referred to as the outsourcing problem, consists of partitioning the set of  customers into two subsets: committed customers and outsourced customers. To this end, we introduce an $n$-dimensional binary variable $x=(x_c)_{c\in\mathcal{C}}$, where $x_c=1$ ($x_c=0$) indicates that customer $c\in\mathcal{C}$ is committed (outsourced).
The subsets of committed and outsourced customers can be defined based on the decision variable $x$ as $\mathcal{C}^o=\{c\in\mathcal{C}|x_c=0\}$ and $\overline{\mathcal{C}}^o=\mathcal{C}\setminus\mathcal{C}^o$.

The second-level problem, which we refer to as the Vehicle Routing Problem with Stochastic Demands (i.e., \vrpx{}),
consists of routing the LSP's fleet to serve the committed customers. Therefore, the second-level problem is stochastic and dynamic, and its solution is expressed in terms of a policy rather than a fixed routing plan. We let $\pi_r$ denote a routing policy used to serve the committed customers $\overline{\mathcal{C}}^o$. Let $w=\{w^c\}_{c\in\overline{\mathcal{C}}^o}$ be the vector of realized customer demands, where each $w^c$ is a demand realization distributed according to $\Gamma^D_c$. For a given routing policy $\pi_r$, committed customer set $\overline{\mathcal{C}}^o$, vehicle $v\in\mathcal{V}$, and demand realization vector $w$, let $T_{\pi_r}(\overline{\mathcal{C}}^o, v, w)$ denote the duration of the route performed by vehicle $v$. Accordingly, we define $R(\overline{\mathcal{C}}^o)$ as the minimum expected routing cost of serving the committed customers, including both travel and overtime costs. We  formalize the \vrp{} as follows:
\begin{align}
    \min_{x} & ~~~ \Psi(\mathcal{C}^o) + R(\overline{\mathcal{C}}^o), \label{eq:obj20} \\
    \textrm{s.t.:}& \nonumber\\
    & \mathcal{C}^o = \{c\in\mathcal{C}|x_c=0\}, \label{eq:cnst21}\\
    & \overline{\mathcal{C}}^o = \{c\in\mathcal{C}|x_c=1\}, \label{eq:cnst22}\\
    & R(\overline{\mathcal{C}}^o) = \min_{\pi_r} \mathbb{E}_{w} \Biggl[ \sum_{v \in \mathcal{V}} \Bigl[ \min\Bigl(T_{\pi_r}(\overline{\mathcal{C}}^o, v, w), ~L\Bigr) + \phi\max\Bigl(T_{\pi_r}(\overline{\mathcal{C}}^o, v, w) - L, ~0\Bigr)  \Bigr] \Biggr], \label{eq:cnst23} \\
    & x_c \in \{0, 1\},~\forall c\in\mathcal{C}.\label{eq:cnst24}
\end{align}

The challenge of formulation \eqref{eq:obj20}-\eqref{eq:cnst24} is the inner optimization problem \eqref{eq:cnst23}, which hides substantial complexity. Indeed, for any fixed first-level decision $x$, evaluating the objective function requires computing $R(\overline{\mathcal{C}}^o)$, that is, solving the second-level problem.
This problem is itself a stochastic, dynamic, and combinatorial optimization problem: routing decisions must adapt in real time to revealed demands, multiple vehicles must be coordinated, and operational features such as vehicle capacity and overtime penalties must be taken into account.
As a result, accurately evaluating a single candidate outsourcing decision is computationally very demanding. The ability to efficiently solve problem \eqref{eq:cnst23} is a key component of the overall solution method.
We model \vrpx{} as an MDP, as detailed in Section \ref{sec:MDP}. This provides a suitable framework for capturing the sequential and stochastic nature of the problem and  for later developing a reinforcement learning approach to approximate routing costs efficiently for the evaluation of the first-level outsourcing decisions.

\subsection{MDP formulation for the \vrpx{}} \label{sec:MDP}
For a given set of committed customers $\overline{\mathcal{C}}^o$, the \vrpx{} dynamically determines the sequence of customer visits and depot returns for each vehicle $v\in\mathcal{V}$ so as to minimize the expected routing and overtime costs. We formulate this problem as an MDP.
We define a \textit{decision epoch} as a point in time at which a vehicle is ready to depart from its current location, that is, either the depot or a customer location.
At each decision epoch $k$, the system state $s_k$ is defined as
\begin{equation}\label{eq:state_cent}
    s_k = (F^{\mathcal{C}}, F^{\mathcal{V}}, t_k),
\end{equation}
where $F^{\mathcal{C}}=[(l_c, h_c, \bar{d}_c, \hat{d}_c)]_{c\in \overline{\mathcal{C}}^o}$ and $F^{\mathcal{V}}=[(l_v, a_v, q_v, g_v)]_{v\in \mathcal{V}}$ represent the state of customers and vehicles, respectively. 
We note that $F^{\mathcal{C}}$ and $F^{\mathcal{V}}$ are matrices with dimensions of $|\overline{\mathcal{C}}^o|\times 4$ and $|\mathcal{V}|\times 4$, respectively.
The last component $t_k$ denotes the time of decision epoch $k$.
The state of each customer $c\in \overline{\mathcal{C}}^o$ is described by its location $l_c$, availability $h_c$, expected demand $\bar{d}_c$, and unserved demand $\hat{d}_c$. We note that $h_c=0$ if customer $c$ is either currently assigned to a vehicle or it is fully served.
The state of each vehicle $v\in \mathcal{V}$ is described by its destination $l_v$, arrival time at destination $a_v$, current available capacity $q_v$, and a binary variable $g_v$ indicating whether the vehicle is in operation. We note that $g_v=0$ indicates that the operation of vehicle $v$ is terminated and it cannot be used in the current or in future decision epochs.
For vehicles that have already arrived at their destination, $l_v$ indicates their current location.

In this formulation we define the action of each vehicle as its next location to visit, namely, a customer or the depot, \citep[see, for example,][]{Maxwell2010,Chen2019,Dastpak2023}.
In particular, we define the joint action at decision epoch $k$ by the $m$-dimensional vector of $y_k=(y_k^1, ..., y_k^m)$, where $y_k^v$ denotes the action of vehicle $v$. 
Let $\bar{\mathcal{V}}_k=\{v\in\mathcal{V}\mid a_v=t_k \land g_v=1\}$ be the set of vehicles that have arrived at their destinations and are available to take an action at decision epoch $k$. We call $\bar{\mathcal{V}}_k$ the set of active vehicles. Since service is assumed to be instantaneous, the arrival time of a vehicle at a customer coincides with the time at which the vehicle is ready to depart.
We now describe the feasible actions for vehicles $v\in\bar{\mathcal{V}}_k$, that is, the possible values of $y_k^v$.
Let $\overline{\mathcal{C}}^o_k = \{ c\in\overline{\mathcal{C}}^o|h_c = 1\}$ be the set of available customers at the decision epoch $k$.
If vehicle $v$ is located at a customer, it has three possible actions. It may visit another customer $c\in\overline{\mathcal{C}}^o_k$ \textit{directly}, denoted by $\mathtt{y}^D_c$, visit another customer  $c\in\overline{\mathcal{C}}^o_k$ \textit{indirectly}, denoted by $\mathtt{y}^I_c$, that is, by first performing an unloading action at the depot, or select the depot $l_0$ as its next destination.
We note that returning to the depot without obligation, that is, not because of a route failure or because no customer is currently available, is interpreted as terminating the trip of the vehicle.
On the contrary, if the vehicle is at the depot because it was obliged to, it may either choose a customer $c$ to visit directly, that is, select action $\mathtt{y}^D_c$, or terminate its operation by choosing action $l_0$.
When the vehicle is at the depot, indirectly visiting a customer $c$, that is, choosing $\mathtt{y}^I_c$, is not applicable.
The joint action space $A(s_k)$ is then defined as follows:
\begin{align} \label{eq:action_space}
    A(s_k) = & \Biggl\{y_k\in \Bigl\{ \{l_0\}\cup \{\mathtt{y}^D_c, \mathtt{y}^I_c| c\in\overline{\mathcal{C}}^o_k\}\Bigr\}^m:\\  
    & y_k^v = l_v & \forall v \in \mathcal{V}\setminus\bar{\mathcal{V}}_k, \label{eq:act_1}\\
    & y_k^v = l_0 &  \forall \{v\in\mathcal{V}|q_v = 0 \lor \overline{\mathcal{C}}^o_k = \emptyset \}, \label{eq:act_3} \\
    & y_k^v \neq \mathtt{y}^I_c & \forall \{v\in\mathcal{V}|l_v=l_0\},\forall c\in\overline{\mathcal{C}}^o_k, \label{eq:act_4} \\
    & y_k^v \neq l_0 & \forall \{v\in\mathcal{V}|\sum_{v\in\mathcal{V}} g_v=1 \land \overline{\mathcal{C}}^o_k \neq \emptyset \}, \label{eq:act_5} \\
    & y_k^v \neq y_k^{v'} & \forall \{v, v' \in \mathcal{V} | v\neq v' \land y_k^v\neq l_0\}. \label{eq:act_6} \Biggr\},
\end{align}
In the action space $A(s_k)$, Condition \eqref{eq:act_1} obliges unavailable  vehicles at time $t_k$ to continue to their current assigned destination.
Condition \eqref{eq:act_3} enforces a return to the depot when a vehicle has no remaining capacity or when no customer is available to visit.
Condition \eqref{eq:act_4} ensures that indirect service of a customer is not allowed when the vehicle is already at the depot.
Condition \eqref{eq:act_5} implies that the vehicle cannot terminate its trip (travel to the depot without obligation), if customers are still available to be served and no other vehicle is in operation.
Finally, Condition \eqref{eq:act_6} ensures that no two vehicles choose the same destination, unless that destination is the depot.

We now describe the evolution of the system from the initial decision epoch onward.
The state at decision epoch zero, denoted by $s_0$, represents the initial configuration of the system before any routing actions are taken. It is constructed from the committed customer set $\overline{\mathcal{C}}^o$, which is determined by the first-level outsourcing decision. At this stage, customer demands have not yet been realized, all vehicles are located at the depot, and no service has started. The initial state $s_0$ is defined as follows:
\begin{equation}\label{eq:initSate}
s_0=([(l_c, 1, \bar{d}_c, ?)]_{c\in \overline{\mathcal{C}}^o}, [(l_0, t_0, Q, 1)]_{v\in\mathcal{V}}, t_0).
\end{equation}
At this decision epoch, the customer's availability $h_c$ and unserved demand $\hat{d}_c$ are initially set to 1 (i.e., available) and ``?" (i.e., unknown), respectively. If customer $c$ is already fully served or a vehicle is assigned to serve it, then $h_c$ takes the value of 0. 
The unserved demand $\hat{d}_c$ is updated to the actual demand upon the first visit.
The state of the fleet is initialized by setting $l_v, a_v, q_v,$ and $g_v$ to $l_0, t_0, Q,$ and 1, respectively.

In state $s_k$, taking an action $y_k\in A(s_k)$ transitions the system to the next state $s_{k+1}$ according to the state transition function $S^M(s_k, y_k, w_{k+1})$.
Recalling that $w$ is a given realization of the stochastic demands, we denote by $w_k\subset w$ the set of demand realizations corresponding to customers that are visited for the first time at decision epoch $k$.

The transition function is divided into two conceptual steps. The first step updates the state $s_k$ to the post-decision state $s_k^y$, which reflects the status of the system immediately after action $y_k$ is taken but before the realization of the exogenous information $w_{k+1}$. The state is updated as follows:
\begin{equation}
    h_c = 0,~ \forall ~ \{c \in \overline{\mathcal{C}}^o| \exists v\in \bar{\mathcal{V}}_k, y_k^v= \mathtt{y}^D_c \lor y_k^v= \mathtt{y}^I_c\},
    \label{eq:pds1}
\end{equation}
\begin{equation}
    g_v=0, ~\forall \{v\in\bar{\mathcal{V}_k}| y^v_k=l_0 \land q_v > 0 \land \overline{\mathcal{C}}^o_k \neq \emptyset \},
    \label{eq:pds2}
\end{equation}
\begin{equation}
    l_v= 
    \begin{cases}
        l_c & y^v_k = \mathtt{y}^D_c \lor \mathtt{y}^I_c, \\
        l_0 & \textrm{otherwise}
    \end{cases},
    \forall ~ v \in \bar{\mathcal{V}}_k,
    \label{eq:pds3}
\end{equation}
\begin{equation}
    a_v = t_k + 
    \begin{cases}
        \tau_{v,0} & y^v_k = l_0 \\
        \tau_{v,c} & y^v_k = \mathtt{y}^D_c \\
        \tau_{v,0} + \tau_{0,c} & y^v_k = \mathtt{y}^I_c
    \end{cases},~ \forall ~ v \in \bar{\mathcal{V}}_k,
    \label{eq:pds4}
\end{equation}
\begin{equation}
    t_{k+1}=\min_{v\in\mathcal{V}} a_v.
    \label{eq:pds5}
\end{equation}
Equation~\eqref{eq:pds1} updates the selected customers as unavailable.
Equation~\eqref{eq:pds2} updates the in-operation status of vehicles that terminate their trip.
Equations~\eqref{eq:pds3} and \eqref{eq:pds4} update the location and arrival time of each vehicle according to its action.
Lastly, Equation~\eqref{eq:pds5} advances the time to the next decision epoch, at which the exogenous information $w_{k+1}$ is revealed.
All other state components remain unchanged.

The second step of the transition function occurs at the beginning of decision epoch $k+1$ and transforms the post-decision state $s_k^y$ into $s_{k+1}$.
Let $\tilde{\mathcal{C}}_{k+1} = \{c\in\overline{\mathcal{C}}^o| \exists v\in\bar{\mathcal{V}}_{k+1} \land l_c=l_v  \}$ be the set of customers that are served at decision epoch $k+1$.
If the actual demand of a customer in $\tilde{\mathcal{C}}_{k+1}$ has not yet been realized (i.e., $\hat{d}_c=?$), then $\hat{d}_c$ is set to the observed demand $w^c_{k+1}$.
Thus,
\begin{equation}
    \hat{d}_c = w^c_{k+1}, ~ \forall ~ \{c\in \tilde{\mathcal{C}}_{k+1}| \hat{d}_c=?\}.
    \label{eq:trns1}
\end{equation}
Let $\eta_v^{k+1}$ denote the demand volume served by vehicle $v\in\bar{\mathcal{V}}_{k+1}$ at decision epoch $k+1$, defined as:
\begin{equation}
    \eta_v^{k+1} = \min \{\hat{d}_{c_v}, q_v\}, ~ \forall ~ \{v \in\bar{\mathcal{V}}_{k+1} | l_v \neq l_0\},
    \label{eq:trns2}
\end{equation}
where $c_v$ denotes the customer $c\in\tilde{\mathcal{C}}_{k+1}$ served by vehicle $v$ (i.e., $l_c=l_v$).
The unserved demands of customers in $\tilde{\mathcal{C}}_{k+1}$ and the available capacity of vehicles $v\in\bar{\mathcal{V}}_{k+1}$ are then updated as follows: 
\begin{equation}
    \hat{d}_{c_v} =\hat{d}_{c_v} - \eta_v^{k+1}, \hspace{5pt} \forall \hspace{3pt} \{v\in\bar{\mathcal{V}}_{k+1} | l_v\neq l_0\}.
    \label{eq:trns3}
\end{equation}
\begin{equation}
    q_v = \begin{cases}
    q_v - \eta_v^{k+1} &  y^v_k = \mathtt{y}^D_c, \\
    Q - \eta_v^{k+1} &  y^v_k = \mathtt{y}^I_c, \\
    Q &  \textrm{otherwise  } .
    \end{cases} \hspace{5pt} \forall \hspace{3pt} v \in \bar{\mathcal{V}}_{k+1}.
    \label{eq:trns4}
\end{equation}
Finally, customers in $\tilde{\mathcal{C}}_{k+1}$ whose demands are not fully served become available again for future service.
Thus,
\begin{equation}
    h_c=1, \hspace{5pt}  \forall \hspace{3pt} \{c \in \tilde{\mathcal{C}}_{k+1}|\hat{d}_c > 0\}.
    \label{eq:trns5}
\end{equation}

We define the cost function as $C(s_k, y_k) = \sum_{v\in \bar{\mathcal{V}_k}} C_v(s_k, y^v_k)$, representing the routing cost incurred by taking action $y_k$ in state $s_k$.
Here, $C_v(s_k, y_k^v)$ denotes the individual routing cost associated with vehicle $v$ taking action $y_k^v$ in state $s_k$, and is defined as follows:
$$C_v(s_k, y^v_k) = \Bigl(\min\{a_v, L\} - \min\{t_k, L\}\Bigr) + \phi \Bigl(\max\{a_v, L\} - \max\{t_k, L\}\Bigr)$$

In this MDP formulation, the routing policy $\pi_r$ determines the action $y_k\in A(s_k)$ when state $s_k$ is observed (i.e., $y_k=\pi_r(s_k)$).
Accordingly, we define the value of being at state $s_k$ as follows:
\begin{equation}\label{eq:value_0}
    V^{\pi_r}(s_k) = \mathbb{E}_{w}\!\left[ C(s_k, y_k) + \gamma V^{\pi_r}(s_{k+1}) \mid s_k \right], \quad \forall s_k \in S_{\overline{\mathcal{C}}^o},
\end{equation}
where $\gamma$ is a discount factor that reflects the relative importance of near-future costs.
In addition, $S_{\overline{\mathcal{C}}^o}$ denotes the state space associated with the committed customer set $\overline{\mathcal{C}}^o$.
We seek an optimal routing policy $\pi_r^*$ that minimizes the value function $V^{\pi_r}(s_k)$:
\begin{equation}\label{eq:rp0}
    \pi_r^* = \arg\min_{\pi_r} V^{\pi_r}(s_k),~\forall s_k\in S_{\overline{\mathcal{C}}^o}.
\end{equation}
The value function corresponding to the optimal policy $\pi_r^*$ is denoted $V^*(s_k)$ for all states $s_k$. For a given committed customer set $\overline{\mathcal{C}}^o$, the value function in Equation \eqref{eq:value_0} and its corresponding decision policy in Equation \eqref{eq:rp0} address the routing problem in Equation \eqref{eq:cnst23}. In particular, $V^*(s_0)$ is used as an approximation of the expected routing cost $R(\overline{\mathcal{C}}^o)$.
The proposed MDP formulation remains computationally intractable because of the size of the state and action spaces. In particular, the associated state-action space is too large to allow an explicit evaluation of the value function $V^{\pi_r}(s_k)$ for every state $s_k\in S_{\overline{\mathcal{C}}^o}$. This challenge is commonly referred to as the curse of dimensionality~\citep{Powerll2022}.

A further difficulty arises in the overall two-level problem, since the committed customer set $\overline{\mathcal{C}}^o$ depends on the first-level outsourcing decision. As a result, a routing policy that is optimal for one committed customer set may no longer be optimal for another, or may not even be defined. Since a given initial customer set $\mathcal{C}$ may induce  an exponential number of committed customer subsets, solving a distinct second-level problem for each such subset is not practical. This is also impractical  when using an iterative algorithm that alternates between  the two-levels of the \vrp{}. Therefore, in Section \ref{sec:generalizedMDP}, we generalize the MDP formulation so as to handle arbitrary committed customer sets. Then, in Section \ref{sec:rl}, we address the curse of dimensionality by developing a deep reinforcement learning algorithm. 

We then focus on the solution of the second-level problem, the \vrpx{}, in Section \ref{sec:solvingVRPX}, and present our ILS algorithm in Section \ref{sec:ILS}.

\section{Solution Method}\label{sec:method}
Our solution method is designed to be executed on each day, before the beginning of the operations. 
The ILS, presented in Section \ref{sec:ILS}, explores and selects candidate outsourcing decisions. However, this algorithm relies on the existence of an \emph{oracle} able to provide an estimation of the expected cost of candidate outsourcing decisions, i.e., the value of the second-level problem. The method to obtain such an oracle is described in Section \ref{sec:solvingVRPX}. First, in Section \ref{sec:generalizedMDP} we generalize the formulation provided in Section~\ref{sec:MDP} so that it operates over a distribution of possible customer subsets, rather than on a specific subset of committed customers $\overline{\mathcal{C}}^o$.
 Thus, eliminating the need to retrain our offline algorithm for each instance. The overarching assumption is that, while the set of customers varies at each day, LSPs typically possess historical information on previous days of operations, which allows them to deduce empirical distributions of customer locations and demands.
Although formally simple, this generalization is conceptually important because it decouples the second-level model from the first-level decision, opening the door to developing a learning-based algorithm able to return a reusable value function that is entirely computed offline.
Second, following \cite{Dastpak2023}, Section \ref{sec:rl} simplifies the obtained formulation by applying the concepts of \emph{consecutive action space} and \emph{observation function}. As it will be clear later on, the resulting simplified MDP formulation allows only one vehicle to be active at a given decision epoch, and adopts an observation function mapping the original large state space into a simplified state tailored to the vehicle that is active at that decision epoch.

Differently from \cite{Dastpak2023} where the observation function is comprised of hand selected features, in Section \ref{sec:GAN} we propose presenting the observation function via a GAT. We then embed the GAT architecture into a deep reinforcement learning approach based on a DQN (Section \ref{sec:QN}). The full network is trained offline, and enables fast and scalable estimation of the routing cost for any subset of committed customer, thus providing the desired oracle to be used within the ILS algorithm.
The details of the ILS are presented in Section \ref{sec:ILS}. In particular, Section \ref{sec:vns} introduces a core implementation of the algorithm while Section \ref{sec:ilsenhancement} introduces an online fine-tuning step of the DQN to enhance the accuracy of the offline-trained value function.

\subsection{Solving the vehicle routing problem with stochastic demands}\label{sec:solvingVRPX}
\subsubsection{A generalized MDP formulation}\label{sec:generalizedMDP}
The MDP formulation presented in Section \ref{sec:MDP} is formalized on a specific customer set $\overline{\mathcal{C}}^o$ with corresponding state space $S_{\overline{\mathcal{C}}^o}$. 
It is important to note that the corresponding value function $V^*(s_k)$ becomes sub-optimal or even undefined for any change in the set of the committed customers $\overline{\mathcal{C}}^o$, which directly impacts $S_{\overline{\mathcal{C}}^o}$.
In our solution method this is an important limitation due to two main factors. First, given a set of customers $\mathcal{C}$ for a given day, we aim at using an ILS to search into possible outsourcing decisions, potentially including any possible subset of $ \mathcal{C}$.  
Solving a separate instance of the problem described in \ref{sec:MDP} for every possible subset of $\overline{\mathcal{C}}^o \subset \mathcal{C}$ is not practical, especially in presence of a limited computational budget, as in our application. This points towards a solution strategy where part of the computations is done offline, before the set of customers $\mathcal{C}$ is known,   which  brings us to the second challenge.
Given that the customer set $\mathcal{C}$ is unknown when performing the offline computations, we need to consider $\mathcal{C}$ as a random variable with a related probability distribution.

The proposed modeling strategy is to generalize the MDP in Section \ref{sec:MDP} by allowing the initial set of committed customers to be a random variable. We formalize this by assuming that the set of daily customers is a random variable distributed as $\mathcal{C} \sim \Gamma^{\mathcal{C}}$. The ILS will then explore the subset of committed customers according to an endogenous conditional probability distribution $\mathcal{\overline{C}}^o~|~\mathcal{C} \sim \Gamma_{ILS}^{\mathcal{\overline{C}}^o|\mathcal{C}}$.
We assume then that
the set of committed customers is distributed according to the marginal distribution $\mathcal{\overline{C}}^o \sim
\Gamma^{\overline{\mathcal C}^{\,o}}
=
\mathbb E_{\mathcal C \sim \Gamma^{\mathcal C}}
\Bigl[
\Gamma_{ILS}^{\overline{\mathcal C}^{\,o}\mid \mathcal C}(\cdot \mid \mathcal C)
\Bigr]$.
In consequence, we rewrite the value function and the optimal routing policy in Equations \eqref{eq:value_0} and \eqref{eq:rp0} as:
\begin{equation}\label{eq:value_1}
    V^{\pi_r}(s_k) = \mathbb{E}_{w} [C(s_k, y_k) + \gamma V^{\pi_r}(s_{k+1})~|~s_k], ~\forall s_k \in \bar{S},
\end{equation}
\begin{equation}\label{eq:rp1}
    \pi_r^* = \arg\min_{\pi_r} V^{\pi_r}(s_k),~\forall s_k\in \bar{S},
\end{equation}
where $\bar{S}$ is the union of state-spaces corresponding to all customer subsets $\overline{\mathcal{C}}^o$ considered under the distribution $\Gamma^{\mathcal{\overline C}^o}$, i.e.,
$\bar{S}=\bigcup_{\overline{\mathcal{C}}^o \sim \Gamma^{\mathcal{\overline C}^o}} S_{\overline{\mathcal{C}}^o}$.

Formulation \eqref{eq:value_1} and \eqref{eq:rp1} is independent of the specific day of operation and the corresponding initial set of customers  $\mathcal{C}$. This fact comes with the additional benefit that once optimized, the optimal value function \eqref{eq:value_1} and policy \eqref{eq:rp1} can be applied to any given day, for any given subset of customers without the need for additional computation. On the other hand, challenges  may come from the estimation of the conditional probability distribution $ \Gamma^{\overline{\mathcal C}^{\,o}}$, given that it is hard to estimate a priori the behavior of the ILS algorithm. In Section \ref{sec:protocol} we explain the choices made in our implementation.

\subsubsection{The consecutive action rule and the observation function \texorpdfstring{\citep{Dastpak2023}}{}}\label{sec:rl}

Solving problem \eqref{eq:value_1} and \eqref{eq:rp1} is extremely challenging, mainly due to the so called curse of dimensionality \citep{Powerll2022}. In particular we face 1) an extremely large action space, due to the potential presence of multiple active vehicles at a given decision epoch, resulting in mutually interdependent actions, see Equations \eqref{eq:action_space} -- \eqref{eq:act_6} describing the mutual decision space, and 2) an extremely large state space, as defined in \eqref{eq:state_cent}.

For a similar problem,  \cite{Dastpak2023} proposed to approximate the original problem formulation by a simpler one obtained by following two main strategies: 1) to restrict the set of feasible actions by imposing a \emph{consecutive action} rule, i.e., by imposing that only one vehicle at time is allowed to be active, and 2) to introduce an aggregation of the original system state called the \emph{observation function}, whose purpose is to filter the information in the original state that is more relevant when establishing the next action to assign to the currently active vehicle. In this section, we briefly describe the application of these ideas to our problem.

The consecutive action rule requires that only one vehicle can be active. This is aligned with the specific characteristic of our problem where it is unlikely that several vehicles simultaneously arrive at their respective destinations and thus are simultaneously active. However, in cases where this happens, such as at decision epoch zero, when a decision must be taken for all vehicles, a random order is imposed among them by assigning each vehicle to a different decision epoch. This results in a singleton active vehicle set $\bar{\mathcal{V}}_k$, with $|\bar{\mathcal{V}}_k|=1$ and $\bar{v} \in \bar{\mathcal{V}}_k$ represents the \textit{active vehicle}. Given the consecutive action rule, the active vehicle $\bar{v}$ for a given state $s_k$ is uniquely defined by the state itself. With this formalism, we set $y_k=(y^{\bar{v}}_k)$ and we define the action space $\bar{A}(s_k)$ as follows:
\begin{align} \label{eq:action_space_v}
    \bar{A}(s_k) = & \Biggl\{y_k^{\bar{v}}\in \Bigl\{ \{l_0\}\cup \{\mathtt{y}^D_c, \mathtt{y}^I_c| c\in\overline{\mathcal{C}}^o_k\}\Bigr\}:\\    
    & y_k^{\bar{v}} = l_0 &  q_{\bar{v}} = 0 \lor \overline{\mathcal{C}}^o_k = \emptyset, \label{eq:act_v_1} \\
    & y_k^{\bar{v}} \neq \mathtt{y}^I_c & l_{\bar{v}}=l_0,\forall c\in\overline{\mathcal{C}}^o_k, \label{eq:act_v_2} \\
    & y_k^{\bar{v}} \neq l_0 & \sum_{v\in\mathcal{V}} g_v=1 \land \overline{\mathcal{C}}^o_k \neq \emptyset, \label{eq:act_v_3} \Biggr\},
\end{align}
where Condition \eqref{eq:act_v_1} forces the active vehicle to visit the depot if either it has no capacity or there is no available customer to serve. Condition \eqref{eq:act_v_2} ensures that indirect visits of customers $c\in\overline{\mathcal{C}}^o_k$ are not feasible if the active vehicle is located at the depot. The last condition implies that the vehicle $\bar{v}$ cannot terminate its trip if it is the only vehicle in-operation and there are still customers to be served.

The second strategy consists in replacing the full system state $s_k$ by a fixed-size aggregated representation, called observation, tailored to the active vehicle $\bar{v}$. Specifically, the observation $o_{k,\bar{v}}$ is obtained via an observation function $O(s_k, \bar{v})$; i.e., $o_{k,\bar{v}} = O(s_k, \bar{v})$. Decisions are then made based on this reduced representation, which yields an approximate MDP formulation, called MDP-CO, on the observation space.

Equations \eqref{eq:pds1}--\eqref{eq:trns5} remain applicable to MDP-CO, if we account for the new definition of $\bar{\mathcal{V}}$. We can then rewrite  equations \eqref{eq:value_1} and \eqref{eq:rp1} by adopting the observation as the representation of the system state:
\begin{equation}
    V^{\pi_r}(o_{k,\bar{v}}) = \mathbb{E}_{w} [C(s_k, y_k) + \gamma V^{\pi_r}(o_{k+1,\bar{v}'})~|~o_{k,\bar{v}}],
    \label{eq:value_mdpco}
\end{equation}
\begin{equation}\label{eq:pi_mdpco}
    \pi_r^* = \arg\min_{\pi_r} V^{\pi_r}(o_{k,\bar{v}}),~\forall o_{k,\bar{v}}\in \Omega,
\end{equation}
where $\bar{v}'$ is the active vehicle in decision epoch $k+1$, and $y_k=\pi_r(o_{k,\bar{v}})$. In addition, $\Omega$ represents the set of all possible observations $o_{k,\bar{v}}$. 
It is worth noting that in equations \eqref{eq:value_mdpco} and \eqref{eq:pi_mdpco}, the value function $V^{\pi_r}(\cdot)$ is defined on the observation $o_{k,\bar{v}}$, whereas in equations \eqref{eq:value_1} and \eqref{eq:rp1}, the same notation is used for the value function defined on the system state $s_k$. We retain this notation for simplicity, since the intended meaning is clear from the context. With this convention, the optimal state-based value $V^*(s_k)$ is approximated by the observation-based value $V^*(o_{k,\bar{v}})$.
The approximated expected routing costs $R(\overline{\mathcal{C}}^o)$ for any tentative set of committed customers $\overline{\mathcal{C}}^o$ is simply retrieved in the value $V^{*} (s_0)$.

\subsubsection{Learning the observation function: the Graph Attention Network}\label{sec:GAN}
The careful design of the observation function is crucial when solving problem \eqref{eq:value_mdpco} and \eqref{eq:pi_mdpco}. \cite{Dastpak2023} proposed an observation function based on hand-crafted features of the full state space, to return a fixed-size vector representation. 
In particular, they replaced the customer component of the system state $F^{\mathcal{C}}$, by two simpler vectors.
The first vector provides a general overview of all customers' state by aggregating information using a grid-like discretization technique. The second vector is a subset of $F^{\mathcal{C}}$ built by filtering only a fixed number of \textit{target customers} chosen according to a heuristic metric. While simple to implement, this observation function has several limitations.
First, since target customers are heuristically selected, good candidate customers may be excluded.  
Second, there is a trade-off in the grid-like discretization of customer information: a finer grid yields higher precision, but also increases the size of the state representation and introduces an additional design parameter, namely the grid resolution.
Third, the observation function $o_{k,\bar{v}}$ is defined relative to the index of the active vehicle $\bar{v}$. Since the vehicles are homogeneous, this representation does not explicitly exploit the permutation symmetry of vehicle identities. As a result, states that are operationally equivalent up to a relabeling of vehicles may still need to be treated separately by the value function approximation. To address these limitations, in this section we propose to replace the hand-crafted observation function in \cite{Dastpak2023}
with a learned representation based on a GAT. A GAT is a neural architecture designed to process graph-structured data by learning how to aggregate information across related entities. This is particularly appropriate in our setting, where decisions depend on both customer and vehicle attributes and on network information such as travel times and spatial proximity.
The attention mechanism allows the model to identify the entities that are most relevant to the active vehicle and to combine their information into a fixed-size observation vector. In this way, the representation is learned directly from data, rather than imposed through manual target-customer selection, grid-based discretization, or an index-dependent treatment of homogeneous vehicles.

The proposed GAT is based on Transformer Encoders introduced in \cite{Vaswani2017}, which use an attention mechanism to generate contextualized representations of their inputs.
In essence, the attention mechanism is a neural operator that takes as input a \textit{query} vector and a \textit{context} matrix, and computes \textit{alignment scores} between the query vector and each row vector of the context matrix.  
The scores are then used to produce a new \textit{context-aware} representation of the query vector, often referred to as the \textit{attention output}. Higher alignment scores correspond to
context rows that receive greater importance in the construction of the query representation.
The resulting attention output is a compact, continuous vector representation of the query, capturing its relationship to the information contained in the context matrix.

In our setting, for example, given a customer and its state representation (the query vector), we may aim to produce a compact alternative representation of that customer relative to the states of all vehicles (the context matrix).
Generally, we may wish to compute representations for all customers (a query matrix) relative to all vehicle states.
More formally, we consider a query matrix  $F \in \mathbb{R}^{|F| \times e_F}$, where each row $ F_i \in \mathbb{R}^{e_F}$ represents a query vector. We seek to compute alternative representations for these query vectors.
To this scope, we consider a context matrix $\bar{F} \in \mathbb{R}^{|\bar{F}| \times e_{\bar{F}}}$, whose rows $\bar{F}_i \in \mathbb{R}^{e_{\bar{F}}}$ are context vectors.
Given a query matrix $F$ and a context matrix $\bar{F}$, the attention mechanism computes a new output matrix $H = \text{Att}(F, \bar{F}) \in \mathbb{R}^{|F| \times e}$, where $e$ is a parameter defining the dimension of the attention output, and each row $H_i$ is the contextualized representation of the corresponding query vector $F_i$.
When $F = \bar{F}$, the mechanism becomes a self-attention operation, where each entity is contextualized relative to all others in the same set. In our setting, for example, we may be interested in obtaining a representation of the customer set relative to the customer set itself.

\begin{figure}[!htbp]
    \centering
    \begin{tikzpicture}
[
    x=0.7cm,
    y=0.7cm,
    font=\footnotesize,
    circlebox/.style={
        draw,
        circle,
        minimum size=5.5mm,
        inner sep=0pt,
        align=center
    },
    dashedbox/.style={
        draw,
        dashed,
        rounded corners,
        inner sep=1pt
    },
    rectbox/.style={
        draw,
        minimum width=1.2cm,
        minimum height=0.45cm,
        inner sep=1pt,
        align=center
    },
    arrow/.style={
        -{Stealth[scale=0.8]},
        thick
    },
    bluearrow/.style={arrow, draw=blue},
    greenarrow/.style={arrow, draw=green!50!black},
    purplearrow/.style={arrow, draw=purple},
    dotstyle/.style={
        font=\bfseries
    }
]

\node[circlebox] (F0) at (-5, 6) {$F_0$};
\node[dotstyle] (dots1) at (-3.5, 6) {$\cdots$};
\node[circlebox] (FF) at (-2, 6) {$F_{|F|}$};
\node[dashedbox, fit=(F0) (dots1) (FF)] {};

\node[circlebox] (Fb0) at (0, 6) {$\bar{F}_0$};
\node[dotstyle] (dots2) at (1.5, 6) {$\cdots$};
\node[circlebox] (FFB) at (3, 6) {$\bar{F}_{\bar{|F|}}$};
\node[dashedbox, fit=(Fb0) (dots2) (FFB)] {};

\node[circlebox, fill=blue!20] (u00) at (-5, 4.5) {$U_0^q$};
\node[circlebox, fill=blue!20] (u0F) at (-2, 4.5) {$U^q_{|F|}$};
\node[circlebox, fill=green!20] (ub10) at (-0.65, 4.5) {$U_0^k$};
\node[circlebox, fill=purple!20] (ub20) at (0.65, 4.5) {$U_0^v$};
\node[circlebox, fill=green!20] (ub1F) at (2.35, 4.5) {$U_{|\bar{F}|}^k$};
\node[circlebox, fill=purple!20] (ub2F) at (3.65, 4.5) {$U_{|\bar{F}|}^v$};

\node[circlebox] (a0) at (-5, 2) {$M_0$};
\node[circlebox] (aF) at (-2, 2) {$M_{|F|}$};

\node[rectbox] (dot0) at (-5, -0.5) {Dot-product};
\node[rectbox] (dotF) at (-2, -0.5) {Dot-product};

\node[circlebox] (H0) at (-5, -2) {$H_0$};
\node[circlebox] (HF) at (-2, -2) {$H_{|F|}$};
\node[dashedbox, fit=(H0) (HF)] {};

\node[dotstyle] at (-3.5, 2) {$\cdots$};
\node[dotstyle] at (-3.5, 4.5) {$\cdots$};
\node[dotstyle] at (1.5, 4.5) {$\cdots$};
\node[dotstyle] at (-3.5, -2) {$\cdots$};

\draw[arrow] (F0) -- (u00);
\draw[arrow] (FF) -- (u0F);
\draw[arrow] (Fb0) -- (ub10);
\draw[arrow] (Fb0) -- (ub20);
\draw[arrow] (FFB) -- (ub1F);
\draw[arrow] (FFB) -- (ub2F);

\draw[bluearrow] (u00) to[out=-90, in=90] (a0);
\draw[bluearrow] (u0F) to[out=-90, in=90] (aF);

\draw[greenarrow] (ub10) to[out=-90, in=40] (a0);
\draw[greenarrow] (ub10) to[out=-90, in=40] (aF);

\draw[purplearrow] (ub20) to[out=-90, in=30] (dot0);
\draw[purplearrow] (ub20) to[out=-90, in=70] (dotF);

\draw[greenarrow] (ub1F) to[out=-140, in=40] (a0);
\draw[greenarrow] (ub1F) to[out=-140, in=40] (aF);

\draw[purplearrow] (ub2F) to[out=-90, in=30] (dot0);
\draw[purplearrow] (ub2F) to[out=-90, in=70] (dotF);

\draw[arrow] (a0) -- (dot0);
\draw[arrow] (aF) -- (dotF);
\draw[arrow] (dot0) -- (H0);
\draw[arrow] (dotF) -- (HF);

\end{tikzpicture}

    \caption{Attention mechanism $Att(F, \bar{F})$}
    \label{fig:attention_gen}
\end{figure}
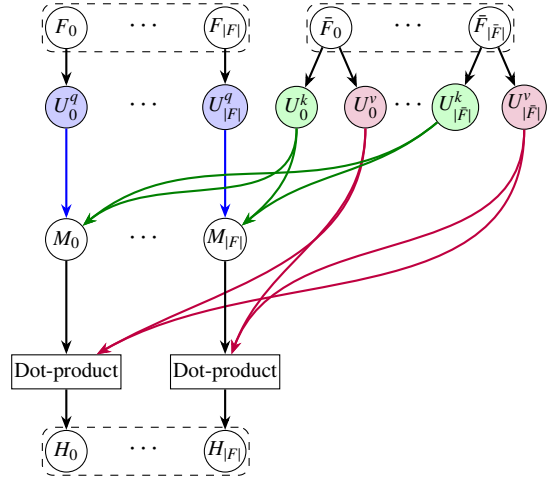

The attention mechanism function $Att(F, \bar{F})$ works as follows: it first projects the query matrix $F$ into a new matrix 
$U^q\in \mathbb{R}^{|F| \times e}$, with $ U^q= F W^q$, and $W^q \in \mathbb{R}^{e_F \times e}$ a learnable weight matrix. Similarly, the context matrix $\bar{F}$ is projected into two separate matrices, $U^k,~ U^v \in \mathbb{R}^{|\bar{F}| \times e}$, with $U^k = \bar{F} W^k$ (keys),  $U^v = \bar{F} W^v$ (values), and $W^k,~W^v \in \mathbb{R}^{e_{\bar{F}} \times e}$ learnable weight matrices.
The alignment between each query vector in $U^q$ and the keys in $U^k$ is computed using scaled dot-product similarity. Specifically, the attention score matrix $M \in \mathbb{R}^{|F| \times |\bar{F}|}$ is calculated as:

\[
M = \text{softmax}\left( \frac{U^q (U^k)^\top}{\sqrt{e}} \right).
\]
The final attention output matrix $H \in \mathbb{R}^{|F| \times e}$ is then computed as
$H = M U^v$,
where each row $H_i$ is a weighted combination of the value vectors in $U^v$, contextualized by the alignment scores of $F_i$ with the components in $\bar{F}$. This completes the computation of $H = \text{Att}(F, \bar{F})$.
Figure \ref{fig:attention_gen} illustrates the schematic structure of the attention mechanism.

Matrices $W^i,\forall i\in\{q,k,v\}$ are the only trainable parameters in the attention mechanism. Their dimensions depend solely on the embedding sizes $e, e_F,$ and $e_{\bar{F}}$, which are fixed. 
Consequently, the number of trainable parameters is independent of $|F|$ and $|\bar{F}|$.
This implies that, once trained, the attention mechanism can process inputs with different sizes (e.g., different number of customers or vehicles) while producing a fixed-size representation of the system state.

\begin{figure}[!ht]
    \centering
    \resizebox{0.98\textwidth}{!}{
    \begin{tikzpicture}[font=\footnotesize]
        \node[anchor=south west, inner sep=0] (image) at (0,0) {\includegraphics[width=0.95\textwidth]{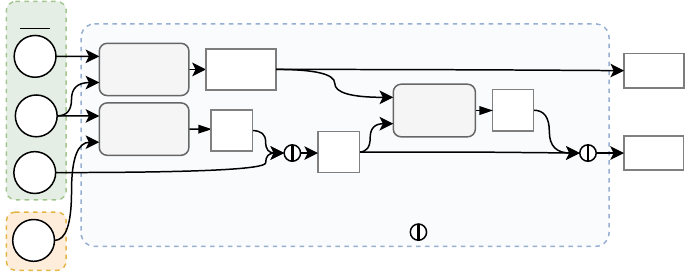}};
        \node[align=center, text width=5cm] at (0.85,5.8) {$s_k$};
        \node[align=center, text width=5cm] at (0.85,5.) {$F^{\mathcal{N}}$};
        \node[align=center, text width=5cm] at (0.85,3.65) {$F^{\mathcal{V}}$};
        \node[align=center, text width=5cm] at (0.85,2.4) {$t_k$};
        \node[align=center, text width=5cm] at (0.8,0.8) {$\bar{v}$};
        \node[align=center, text width=1.5cm] at (3.2,4.65) {Node Embedder};
        \node[align=center, text width=1.5cm] at (3.2,3.3) {Vehicle Embedder};
        \node[align=center, text width=1.5cm] at (9.7,3.8) {Graph Embedder};
        \node[align=center, text width=5cm] at (5.4,4.65) {$H^{\mathcal{NV}}$};
        \node[align=center, text width=5cm] at (5.3,3.3) {$H^{\bar{v}}$};
        \node[align=center, text width=5cm] at (7.6,2.8) {$H^{\bar{v}, t}$};
        \node[align=center, text width=5cm] at (11.5,3.8) {$H^{\mathcal{G}}$};
        \node[align=center, text width=5cm] at (14.8,4.65) {$H^{\mathcal{NV}}$};
        \node[align=center, text width=5cm] at (14.7,2.8) {$o_{k,\bar{v}}$};
        \node[align=center, text width=5cm] at (11.6,1.) {:Concatenates input vectors};
        \node[align=center, text width=5cm] at (11.8,5.9) {\underline{Observation Function}};
    \end{tikzpicture}
    }
    \caption{Schematic structure of proposed Observation Function}
    \label{fig:obsfunc}
\end{figure}
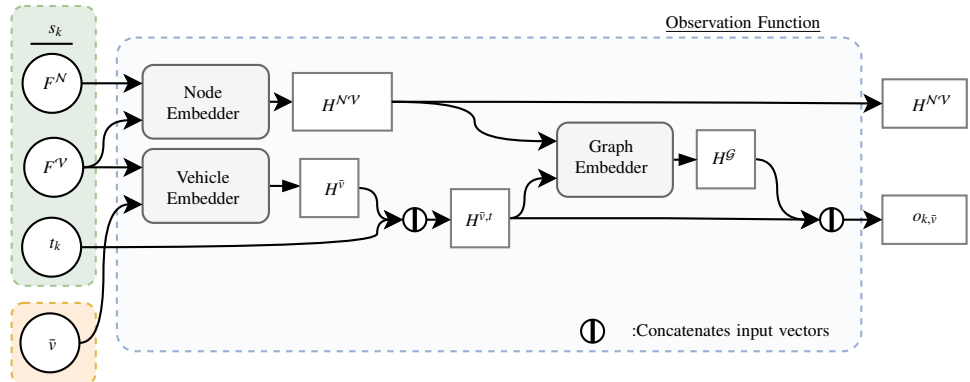

\begin{figure}[!ht]
    \centering
    \subfigure[]{
    \begin{tikzpicture}[
        font=\footnotesize,
        node distance=0.75cm and 1cm,
        every node/.style={},
        round/.style={circle, draw, inner sep=2pt},
        rect/.style={rectangle, draw, inner sep=4pt}
    ]
    
    \node[round] (Fn1) {$F^{\mathcal{N}}$};
    \node[round, right=of Fn1] (Fn2) {$F^{\mathcal{N}}$};
    \node[round, right=of Fn2] (Fn3) {$F^{\mathcal{N}}$};
    \node[round, right=of Fn3] (Fv2) {$F^{\mathcal{V}}$};
    \node[rect, below=of $(Fn1)!0.5!(Fn2)$] (Att1) {$Att(F^{\mathcal{N}}, F^{\mathcal{N}})$};
    \node[rect, below=of $(Fn3)!0.5!(Fv2)$] (Att2) {$Att(F^{\mathcal{N}}, F^{\mathcal{V}})$};
    \node[round, below=of Att1] (Hn) {$H^{\mathcal{N}}$};
    \node[round, below=of Att2] (Hv) {$H^{\mathcal{V}}$};
    \node[draw, circle, minimum size=0.7cm, below=0.5cm of $(Hn)!0.5!(Hv)$] (plus) {$|$};
    \node[round, below=of plus] (Hnv) {$H^{\mathcal{N}\mathcal{V}}$};
    
    \draw[-{Stealth}] (Fn1) -- (Att1);
    \draw[-{Stealth}] (Fn2) -- (Att1);    
    \draw[-{Stealth}] (Fn3) -- (Att2);
    \draw[-{Stealth}] (Fv2) -- (Att2);    
    \draw[-{Stealth}] (Att1) -- (Hn);
    \draw[-{Stealth}] (Att2) -- (Hv);    
    \draw[-{Stealth}] (Hn) -- (plus);
    \draw[-{Stealth}] (Hv) -- (plus);    
    \draw[-{Stealth}] (plus) -- (Hnv);    
    \end{tikzpicture}
    \label{fig:node_emb}}\hfill
    \subfigure[]{
    \begin{tikzpicture}[
        node distance=0.75cm and 1.cm,
        every node/.style={},
        round/.style={circle, draw, inner sep=2pt},
        rect/.style={rectangle, draw, inner sep=4pt}
    ]
    
    \node[round] (Fv_bar) {$F^{\mathcal{V}}_{\overline{v}}$};
    \node[round, right=of Fv_bar] (Fv) {$F^{\mathcal{V}}$};
    \node[rect, below=of $(Fv_bar)!0.5!(Fv)$] (Att) {$Att(F^{\mathcal{V}}_{\overline{v}}, F^{\mathcal{V}})$};
    \node[round, below=of Att] (Hv_bar) {$H^{\bar{v}}$};
    
    \draw[-{Stealth}] (Fv_bar) -- (Att);
    \draw[-{Stealth}] (Fv) -- (Att);
    \draw[-{Stealth}] (Att) -- (Hv_bar);
    
    \end{tikzpicture}
    \label{fig:vehicle_emb}}\hfill
    \subfigure[]{
    \begin{tikzpicture}[
        node distance=0.75cm and 1cm,
        every node/.style={},
        round/.style={circle, draw, inner sep=2pt},
        rect/.style={rectangle, draw, inner sep=4pt}
    ]
    
    \node[round] (Hz) {$H^{\bar{v}, t}$};
    \node[round, right=of Hz] (Hnv) {$H^{\mathcal{NV}}$};
    \node[rect, below=of $(Hz)!0.5!(Hnv)$] (Att) {$Att(H^{\bar{v}, t}, H^{\mathcal{NV}})$};
    \node[round, below=of Att] (Hg) {$H^{\mathcal{G}}$};
    
    \draw[-{Stealth}] (Hz) -- (Att);
    \draw[-{Stealth}] (Hnv) -- (Att);
    \draw[-{Stealth}] (Att) -- (Hg);
    
    \end{tikzpicture}
    \label{fig:graph_emb}}
    \caption{(a) Node Embedder, (b) Vehicle Embedder, (c) Graph Embedder}
    \label{fig:embedders}
\end{figure}
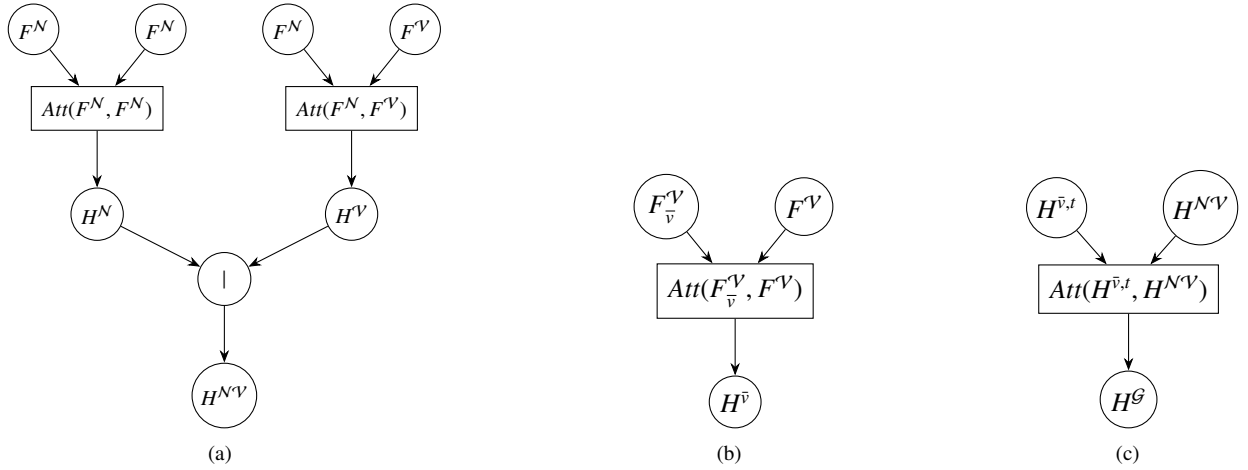

In the following, we present the implementation of the observation function $O(s_k, \bar{v})$ via the  attention mechanism. The general overview is illustrated in Figure \ref{fig:obsfunc}.
The proposed implementation is comprised of three main attention mechanism blocks, namely Node Embedder, Vehicle Embedder, and Graph Embedder. Each of these blocks are illustrated in Figure \ref{fig:embedders}.
The Node Embedder, illustrated in Figure~\ref{fig:node_emb}, takes as input the states of all customers and vehicles. 
For implementation purposes, we assume that the number of customers in a given operational day does not exceed a fixed upper bound, denoted by $n_{\max}$; that is, $n \leq n_{\max}$ for any realization of $\mathcal{C}$. This value is specified in the experimental setup described in Section~\ref{sec:instances}. 
Accordingly, we define a graph
$\mathcal{G}=(\mathcal{N}, \mathcal{E})$, where
$\mathcal{N} = \{0, 1, \ldots, n_{\max}\}$ is the set of $n_{\max} + 1$ nodes, which includes the depot, customers $\overline{\mathcal{C}}^o$, and dummy nodes if $|\overline{\mathcal{C}}^o| < n_{\max}$. Furthermore $\mathcal{E}$ is the set of edges between every two nodes in $\mathcal{N}$, representing potential consecutive visits of a given vehicle. 
We further define $F^{\mathcal{N}} \in \mathbb{R}^{(n_{\max}+1) \times 4}$ as the fixed-size matrix that includes the state representation of the depot and all potential customer nodes at a given decision epoch. 
In this matrix, the first row corresponds to the state of the depot, rows $1$ to $|\overline{\mathcal{C}}^o|$ correspond to $F^{\mathcal{C}}$, and rows $|\overline{\mathcal{C}}^o| + 1$ to $n_{\max}$ are padded with zero vectors (i.e., $(0, 0, 0, 0)$). 
Additionally, a masking strategy is applied to $F^{\mathcal{N}}$ such that the feature vectors of unavailable customers (i.e., $c \in \overline{\mathcal{C}^o} \setminus \overline{\mathcal{C}}^o_k$) are set to zero.
This masking strategy ensures that both unavailable and non-existent nodes are represented by zero vectors, which effectively removes their influence in the attention mechanism. The node embedder also takes
as input the state of all vehicles, represented by $F^{\mathcal{V}} \in \mathbb{R}^{|\mathcal{V}| \times 4}$.

Structurally, the Node Embedder features two parallel attention blocks that are later concatenated into a single attention output.
The first self-attention block, $Att(F^{\mathcal{N}}, F^{\mathcal{N}})$, generates a matrix $H^{\mathcal{N}} \in \mathbb{R}^{(n_{\max}+1)\times e}$, where each row $H^{\mathcal{N}}_i$ is a representation 
of node $i \in \mathcal{N}$ contextualized with respect to  the other customers.
The second attention block, $Att(F^{\mathcal{N}}, F^{\mathcal{V}})$, generates a matrix $H^{\mathcal{V}} \in \mathbb{R}^{(n_{\max}+1)\times e}$, where each row $H^{\mathcal{V}}_i$ is a representation of node $i$ contextualized with respect to the fleet of vehicles.
The final attention output $H^{\mathcal{NV}}$ has dimension $(n_{\max}+1) \times 2e$ and is obtained by concatenating $H^{\mathcal{N}}$ and $H^{\mathcal{V}}$. 
As a result, each node in $\mathcal{N}$ is now represented by a vector of size $1 \times 2e$, providing an enriched representation that captures its relationships with both other nodes in the graph and all vehicles in the system.

The Vehicle Embedder, illustrated in Figure~\ref{fig:vehicle_emb},
takes the state of the active vehicle $F^{\mathcal{V}}_{\bar{v}}$ as the query vector, and state of all vehicles $F^{\mathcal{V}}$ as the context matrix. The resulting attention mechanism $Att(F^{\mathcal{V}}_{\bar{v}},F^{\mathcal{V}})$
outputs a vector $H^{\bar{v}} \in \mathbb{R}^{1\times e}$ representing the active vehicle contextualized with respect to the rest of the fleet. Importantly, this operation does not assume any particular ordering of the vehicles. This ensures permutation invariance, and thus accounts for symmetry among homogeneous vehicles.

The final component of the observation function is the Graph Embedder. As illustrated in Figure~\ref{fig:graph_emb}, it takes as input the active vehicle embedding $ H^{\bar{v}}$, the current time $t_k$, and the node embedding $H^{\mathcal{NV}}$. In particular, the query vector $H^{\bar{v}, t} = [H^{\bar{v}} | t_k]$
concatenates the active vehicle's embedding with the current time. The context matrix is provided by the node embedding  $H^{\mathcal{NV}}$. The corresponding attention mechanism 
$Att(H^{\bar{v}, t}, H^{\mathcal{NV}})$ produces a graph-level embedding $H^{\mathcal{G}} \in \mathbb{R}^{1 \times e}$ providing a representation of the active vehicle contextualized with respect to the embedding of the full set of nodes in the graph. We note that the rationale of including $t_k$ in the query vector, is to allow the attention mechanism to assess the relevance of each node by considering the remaining operational time.

Finally, we define the active vehicle's observation of a given state $s_k$ as:
\begin{equation}
    o_{k,\bar{v}} = [H^{\mathcal{G}}|H^{\bar{v}, t}].
\end{equation}

\subsubsection{The deep Q-learning Algorithm} \label{sec:QN}
In this section, we apply the Q-learning algorithm, which is a model-free version of approximate value iteration~\citep{Watkins1992}, to solve the MDP-CO, specialized to the observation function presented in the previous section. We call the corresponding algorithm \dqnco{}.
In our setting the Q-factor $Q^{\pi_r}(o_{k,\bar{v}},y_k)$ represents the expected cost-to-go of performing the action $y_k$ when $o_{k, \bar{v}}$ is observed, and then applying policy $\pi_r$. We formally define it as follows:
\begin{equation}
    Q^{\pi_r}(o_{k,\bar{v}},y_k) = \mathbb{E}_{w} [C(s_k, y_k) + \gamma V^{\pi_r}(o_{k+1,\bar{v}'})~|~o_{k,\bar{v}},~y_k],
    \label{eq:qvalue_mdpco}
\end{equation}
where $\bar{v}'$ is the active vehicle in decision epoch $k+1$.
We then rewrite the value function in Equation \eqref{eq:value_mdpco} in terms of Q-factors as follows:
\begin{equation} \label{eq:value_q_value}
    V^{\pi_r}(o_{k,\bar{v}}) = Q^{\pi_r}(o_{k,\bar{v}}, \pi_r(o_{k,\bar{v}})), ~\forall o_{k,\bar{v}}\in \Omega.
\end{equation}

\noindent Therefore, the optimal policy $\pi_r^*$ is rewritten as:
\begin{equation} \label{eq:qpolicy}
    \pi_r^* =\arg\min_{\pi_r} Q^{\pi_r}(o_{k,\bar{v}}, \pi_r(o_{k,\bar{v}})), ~\forall o_{k,\bar{v}}\in \Omega.
\end{equation}
Consistently with the convention previously introduced for the value function $V^{\pi_r}(\cdot)$, we use the notation $Q^{\pi_r}(\cdot)$ both for Q-factors defined on the system state $s_k$ and for Q-factors defined on the observation $o_{k,\bar{v}}$.
The Q-factors corresponding to the optimal policy $\pi_r^*$ are denoted $Q^*(\cdot)$.
Thus, the optimal state-based Q-factor $Q^*(s_k,y_k)$ is approximated by the observation-based Q-factor $Q^*(o_{k,\bar{v}},y_k)$.

The optimal policy $\pi^*_r$ under Q-learning is found through a sequence of interactions of an agent with the environment, the so-called experiences, and by iteratively updating an estimate of the Q-factors.
We define an \textit{experience} as a tuple of
$(s_k,~y_k, ~c_k,~s_{k+1})$, describing the agent at the decision epoch $k$ in state $s_k$, taking action $y_k$ and observing a cost $c_k=C(s_k, y_k)$. The agent then transitions to the state $s_{k+1}$.
The proposed \dqnco{} is illustrated in Figure \ref{fig:decdqnoverview}. The Q-network receives as input the observation $o_{k,\bar{v}}$ and the embedding $H^{\mathcal{NV}}_i$ of a node $i \in \mathcal{N}$ representing the next visit, i.e., $y_k$, and returns the corresponding Q-factors.
More specifically, as defined in Equations \eqref{eq:action_space_v}--\eqref{eq:act_v_3}, the active vehicle may travel to a customer $c\in\overline{\mathcal{C}}^o$ either directly, denoted as $\mathtt{y}^D_c$, or indirectly (by passing at the depot), denoted as $\mathtt{y}^I_c$.
Accordingly, we associate two Q-factors with each  customer $c\in\overline{\mathcal{C}}^o$.
Coherently, the proposed Q-network returns two Q-factors ($Q^D_c, Q^I_c$) for each input customer $c$. 
Note that to ensure action feasibility, we set $Q^I_c$ to $+\infty$, if serving that customer indirectly is not feasible, i.e., $\mathtt{y}^I_c\not\in \bar{A}(s_k)$. 
To compute the Q-factor of returning to the depot, we follow a similar procedure, by passing the node embedding of the depot, i.e., $H^{\mathcal{NV}}_0$, to the Q-network. In this case, the second $Q^I_0$ is not applicable, thus we set it to $+\infty$.
\definecolor{softgreen}{RGB}{205,224,178}
\definecolor{softorange}{RGB}{250,224,168}
\definecolor{softblue}{RGB}{211,227,244}
\definecolor{softblueborder}{RGB}{139,170,209}
\definecolor{softpurple}{RGB}{226,210,237}
\definecolor{softpurpleborder}{RGB}{150,110,180}
\definecolor{boxgray}{RGB}{128,128,128}

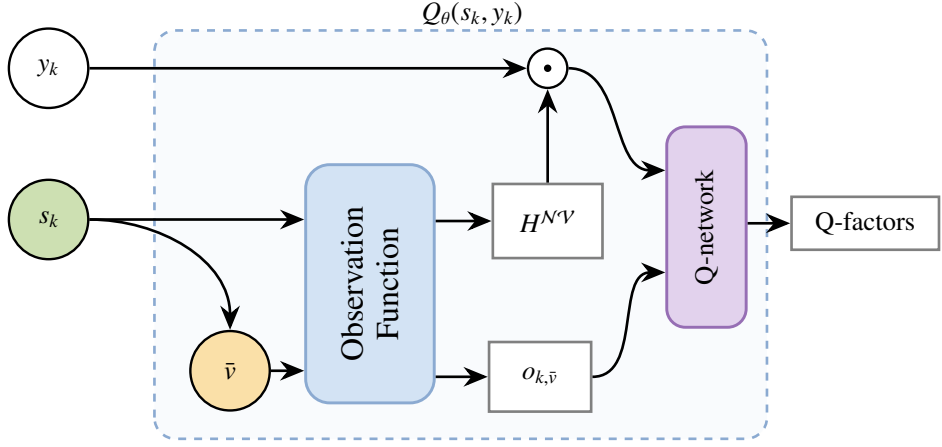
\begin{figure}
    \centering
    \begin{tikzpicture}[
      >={Stealth[length=3.2mm,width=2.6mm]},
      line width=1pt,
      circnode/.style={circle, draw=black, line width=0.9pt, minimum size=1.05cm, inner sep=0pt},
      graybox/.style={draw=boxgray, fill=white, line width=1pt, inner sep=3pt},
    ]

    \draw[softblueborder, dashed, line width=1pt, fill=softblue!45, fill opacity=0.3, rounded corners=8pt]
      (0.2,-0.1) rectangle (8.3,5.3);
    \node[anchor=north east, inner sep=6pt] at (5.3,5.9) {$Q_\theta(s_k,y_k)$};
     
    \node[circnode, fill=white]     (yk) at (-1.2, 4.8) {$y_k$};
    \node[circnode, fill=softgreen] (sk) at (-1.2, 2.8) {$s_k$};
    \node[circnode, fill=softorange](vb) at (1.2, 0.8) {$\bar{v}$};
     
    \node[draw=softblueborder, fill=softblue, line width=1.1pt, rounded corners=9pt,
          minimum width=1.7cm, minimum height=3.15cm] (obs) at (3.05, 1.95) {};
    \node at (obs) {\rotatebox{90}{\large\shortstack[c]{Observation\\[2pt]Function}}};
     
    \node[graybox, minimum width=1.45cm, minimum height=0.98cm] (H) at (5.4, 2.78)
          {$H^{\mathcal{NV}}$};
    \node[graybox, minimum width=1.35cm, minimum height=0.92cm] (o) at (5.3, 0.72)
          {$o_{k,\bar{v}}$};
     
    \node[circle, draw=black, line width=0.9pt, minimum size=0.52cm, inner sep=0pt, fill=white]
          (op) at (5.4, 4.8) {};
    \fill (op) circle (1.4pt);
     
    \node[draw=softpurpleborder, fill=softpurple, line width=1.1pt, rounded corners=7pt,
          minimum width=1.05cm, minimum height=2.65cm] (qnet) at (7.5, 2.7) {};
    \node at (qnet) {\rotatebox{90}{Q-network}};
     
    \node[graybox, minimum width=1.95cm, minimum height=0.7cm] (qval) at (9.6, 2.75)
          {Q-factors};

    \draw[->] (yk.east) -- (op.west);
    \draw[->] (sk.east) -- (obs.west |- sk);
    
    \draw[->] (vb.east) -- (obs.west |- vb);
    \draw[->] (obs.east |- H) -- (H.west);
    \draw[->] (obs.east |- o) -- (o.west);
    \draw[->] (H.north) -- (op.south);
     
    \draw[->] (op.east) to[out=0, in=180] ($(qnet.west)+(0,0.75)$);
    \draw[->] (o.east)  to[out=0, in=180] ($(qnet.west)+(0,-0.6)$);
    \draw[->] (sk.east) to[out=0, in=90] ($(vb.north)+(0,0)$);
     
    \draw[->] (qnet.east |- qval) -- (qval.west);
     
    \end{tikzpicture}
    \caption{Schematic overview of the proposed \dqnco{}}
    \label{fig:decdqnoverview}
\end{figure}

We designed the Q-network as a two-layer artificial neural network, connected by Rectified Linear Unit (ReLU) activation functions, denoted by $\sigma$. Thus, Q-factors for a given observation $o_{k,\bar{v}}$ and a node $i\in\mathcal{N}$ with the corresponding node embedding of $H^{\mathcal{NV}}_i$ can be computed as follows: 
\begin{equation}    [Q^D_i,Q^I_i]=\sigma\Bigl(\sigma\bigl([o_{k,\bar{v}}|H^{\mathcal{NV}}_i]~W^Q_1 + b^Q_1\bigr)~W^Q_2 + b^Q_2\Bigr),
\end{equation}
where $W^Q_1, W^Q_2, b^Q_1,$ and $b^Q_2$ are trainable parameters.

Since both the observation function and the Q-network are based on neural networks, we train both networks together in a single deep neural network structure characterized by a parameter $\theta$.
We thus introduce a function $Q_\theta(s_k,y_k)$ representing the combination of the observation function and the Q-network, where $\theta$ expresses the dependence on the values of the trainable parameters. For a given $\theta$, $Q_\theta(s_k,y_k)$ represents an estimate of the Q-factors.
The training procedure aims at minimizing a loss function.
We define the loss function as the Mean Squared Error (MSE) of $\Delta$, where $\Delta$ is the temporal difference between the approximated Q-factor given by the network $Q_{\theta}(s_k, y_k)$ and the estimated expected Q-factor $\hat{Q}(s_k, y_k)$ computed for each experience $(s_k, y_k, c_k, s_{k+1})$.
$\hat{Q}(s_k, y_k)$ is defined as follows:
\begin{equation}\label{eq:singlevcoaobs}
    \hat{Q}(s_k, y_k) = c_k + \gamma\min_{y_{k+1}\in \bar{A}(s_{k+1})} Q_\theta(s_{k+1}, y_{k+1}).
\end{equation}
Let $\theta^*$ be a minimizer of the loss function $\Delta$. We then approximate the optimal Q-factor $Q^*(s_k,y_k)$ by $Q_{\theta^*}(s_k,y_k)$.

After exploring different enhancement techniques, preliminary experiments revealed that using the Double Q-network~\citep{VanHasselt2016} and replay memory help our method to reach better convergence, as well as improving its overall performance. 
The double Q-network considers a secondary network parameterized by $\bar{\theta}$. This network enables us to evaluate the expected Q-factors ($\hat{Q}$) based on a stable target, resulting in a more stable behavior in the training process. 
Therefore, we replace the secondary network with parameters $\bar{\theta}$, by a copy of the primary network with parameters $\theta$, every $P^{d}$
trials. The traditional Q-network methods train the network for the most recent experiences. However, the memory replay technique buffers new experiences into a fixed-size list, denoted as $B$, and enables us to train the network periodically over a subset of experiences randomly sampled from the buffer list. The buffer list has the FIFO property which discards older experiences when its capacity is reached. With the probability $P^t$ at each decision epoch, we uniformly sample a minibatch of experiences, denoted  $\tilde{B}$, from $B$ and use it to update the network parameters $\theta$. Using the replay memory~\citep{Mnih2015} has a number of benefits, including breaking the correlation between successive experiences and giving more opportunities to those that occur often. 
The choice of the enhancement techniques explained above is determined empirically by performing a series of preliminary experiments.
We detail the \dqnco{} in Algorithm~\ref{alg:DecDQN}.

\begin{algorithm}[!ht] 
 \SetAlgoLined
 \DontPrintSemicolon
 \LinesNumbered
 \SetKwInOut{Input}{input}
 \SetKwInOut{Output}{output}
 Initialize $\theta, \bar{\theta}$ \;
 
 \While{trials $<$ Trials$\_$MAX}
 { 
 Observe a new set of customers $\overline{\mathcal{C}}^o$\;
 Generate a sample demand scenario $w$\;
 $k\leftarrow 0$ \;
 
 \While{true}{
 $\bar{v} \leftarrow$ $\bar{\mathcal{V}}_k[0]$\;
 
 \If{$k > 0$}{
    Given the demand scenario $w$, transition from $s_{k-1}^y$ to $s_k$ \; 
        
    Append $(s_{k-1}, y_{k-1}, c_{k-1}, s_k)$ to $B$\;
 }
 
 $y_k=\begin{cases}
     \textrm{random action in }\bar{A}(s_k) & rand[0,1] < \epsilon \\
     \arg\min_{y\in \bar{A}(s_k)} Q_\theta(s_k, y) & \textrm{otherwise}
     \end{cases}$\;
 $c_k = C(s_k, y_k)$\;
 
 Transition from $s_k$ to $s_k^y$\;

 \If{$rand[0,1] < P^{t}$}{
 Sample the minibatch $\tilde{B}$ from $B$\;

 $\Delta \gets [0]_{|\tilde{B}|}$\;
 $i \gets 0$\;
 \ForEach{$(s, y, c, s')\in \tilde{B}$}{
 \eIf{$s'$ is terminal}{
    $\Delta[i] = Q_\theta(s, y) - c$\;
 }{
    $\Delta[i] = Q_\theta(s, y) - [c + \gamma \min_{y'\in \bar{A}(s')}Q_{\bar{\theta}}(s', y')]$
 }
 $i \gets i+1$\;
 }
 Compute the loss function and update the network $\theta$: $Loss=\mathbb{E}_{\tilde{B}}(\Delta^2), \hspace{5pt} \theta = \theta - \alpha \nabla_{\theta} Loss$\;
 
 }
 \If{$k=K$}{
 break
 }
 Set $k = k + 1$\;
 }
 
 Set $trials = trials + 1$\;
 Every $P^d$ $trials$, set $\bar{\theta} = \theta$\;
 Decay the learning rate $\alpha$ and the exploration rate $\epsilon$\;
 }
 
 $\textbf{return}~ Q_\theta(., .)$\;
 \caption{\dqnco{}}
 \label{alg:DecDQN}
\end{algorithm}

\subsection{The ILS algorithm}\label{sec:ILS}
\subsubsection{The core implementation}\label{sec:vns}
In this section we present a meta-heuristic algorithm to solve the \vrp{}. Our method is based on the principles of the Iterated Local Search. The main idea is that the ILS searches over the set of possible outsourced customer subsets $\mathcal{C}^o$ and takes advantage of the offline learning algorithm, described in Section~\ref{sec:solvingVRPX}, to obtain
nearly-instantaneous routing cost estimates of the corresponding committed customer subset $\overline{\mathcal{C}}^o$. A general overview of the algorithm is provided in Algorithm \ref{alg:ILS}.

The ILS takes the customer set $\mathcal{C}$ as input and produces an outsourcing solution $x^*$ as output. It begins by generating an initial solution $x$ defining a subset of  $\mathcal{C}$. The \texttt{initial\_solution()} function assigns $x_c=1$ to a random subset of customers in $\mathcal{C}$ and $x_c=0$ for the remaining customers.
The $\texttt{evaluate\_solution}$$(x)$ function estimates the objective value of solution $x$ by replacing the expected routing cost $R({\cal\bar{C}}^o)$ in \eqref{eq:obj20} and \eqref{eq:cnst23}, by its approximation $V^*(s_0)$, which is computed offline and retrieved by solving $V^*(s_0) = \min_{y_k\in \bar{A}(s_0)} Q_\theta(s_0, y_k)$, as defined in \eqref{eq:value_q_value}. Specifically, 
$$\texttt{evaluate\_solution}(x) = \Psi(\mathcal{C}^o) + V^*(s_0).$$

The ILS algorithm then employs a Variable Neighborhood Search (VNS) procedure (lines 7--16 in Algorithm \ref{alg:ILS}). 
The VNS algorithm seeks the local optimum around $x$ by exploring two neighborhoods. First, the ``Add'' neighborhood structure generates neighbors of solution $x$ by converting an outsourced customer ($x_c=0$) in $x$ to a committed customer ($x_c\leftarrow 1$). Second, the ``Swap'' neighborhood structure generates neighbors of $x$ by simultaneously converting an outsourced customer ($x_c=0$) to a committed customer ($x_c\leftarrow 1$) and setting a committed customer ($x_{c'}=1$) to an outsourced customer ($x_{c'}\leftarrow 0$). The algorithm begins with the ``Add'' neighborhood (lines 8--11), identifying the best neighborhood solution $(x', f')$, and replacing it with the current solution $(x, f)$ if it has a lower objective value (i.e., $f' < f$). If no improvement occurs, the algorithm explores the 'Swap' neighborhood of the current solution $x$ (lines 12--15) and repeats the procedure. If no improvement is found in both neighborhoods, the VNS loop terminates, and solution $(x, f)$ is returned.
At the end of the VNS step, we devise an optional step, called \emph{online fine-tuning} (lines 17--20), which is detailed in Section \ref{sec:ilsenhancement}.
In the next step (lines 21--23), the current solution $(x, f)$ is compared with the best solution $(x^*, f^*)$. If $f < f^*$, the best solution $x^*$ is updated; otherwise, we increment $no\_change$ by one. The VNS algorithm is repeated while $itr < \xi$ and $no\_change < 15$, each time perturbing the best solution $x^*$ to generate a new solution $x$. The \texttt{perturbation}$(x^*, \kappa)$ function shakes the best solution by randomly removing $\kappa$ customers from $\overline{\mathcal{C}}^o$ and assigning them to the outsourced set. This function reduces the risk of getting stuck in local optima. 
\begin{algorithm}[!ht]
\DontPrintSemicolon
 \SetKwInOut{Input}{input}
 \SetKwInOut{Output}{output} 
  $(x^*, f^*)\leftarrow (\emptyset, \infty)$\;
  $itr\leftarrow 0, no\_change\leftarrow 0$\;
 \While{$itr < \xi$ and $no\_change < 15$}
 {
    $\kappa \leftarrow \min(\lceil \frac{no\_change + 1}{5} \rceil, 3)$\tcp*{$\kappa\in\{1, 2, 3\}$}
    $x\leftarrow \begin{cases}
        \texttt{initial\_solution}(\mathcal{C}) & x^*=\emptyset \\ 
        \texttt{perturbation}(x^*, \kappa) & \textrm{otherwise}
    \end{cases}$\tcp*{generate $x$ by initializing or shaking $x^*$}
    $f \leftarrow \texttt{evaluate\_solution}(x)$\;
    \While{true}
    {
        $(x', f')\leftarrow \texttt{best\_improvement}(x, \textrm{Add})$\tcp*{explore the `Add' neighborhood of $x$}
        \If{$f'<f$}{
            $(x, f)\leftarrow (x', f')$\tcp*{update the current solution}
            \textrm{Continue}\;
        }

        $(x', f')\leftarrow \texttt{best\_improvement}(x, \textrm{Swap})$\tcp*{explore the `Swap' neighborhood of $x$}
        \If{$f'<f$}{
            $(x, f)\leftarrow (x', f')$\tcp*{update the current solution}
            \textrm{Continue}\;
        }
        \textrm{break}\;
    }
    \If{online\_fine\_tuning}
    {
    $f^r\leftarrow \texttt{simulate\_routing}(x)$\tcp*{estimate the routing costs via simulation}
    $f\leftarrow \Psi(\{c\in\mathcal{C}|x_c=0\}) + f^r$\tcp*{update the costs of solution $x$}

    $\texttt{fine\_tuning}(x, f^r)$ \tcp*{fine-tune the value function}
    }
    
    \If{$f<f^*$}{
        $(x^*, f^*)\leftarrow (x, f)$\tcp*{update the best solution}
        $no\_change \leftarrow 0$\;
    }\Else{
    $no\_change \leftarrow no\_change + 1$\;
    }
    $itr \leftarrow itr + 1$\;
 }
 
 \textbf{return} $(x^*, f^*)$\;
 \caption{Iterated Local Search}
 \label{alg:ILS}
\end{algorithm}

\subsubsection{ILS algorithm enhancement}\label{sec:ilsenhancement}
Preliminary experiments revealed some differences between routing costs estimated via the value function $V^*(s_0)$ \eqref{eq:value_q_value} and those computed by simulating the routing policy $\pi_r^*$ \eqref{eq:qpolicy} on the customer set $\overline{\mathcal{C}}^o$. 
Specifically, on small instances the observed gaps were relatively modest on average, however, for some $\overline{\mathcal{C}}^o$, this gap was significant. As a result, some outsourcing decisions looked considerably more attractive than what they really were. We noticed that the ILS algorithm sometimes invested considerable effort exploring neighborhoods of these solutions, disregarding potentially better alternatives.
Lines 17--20 of Algorithm \ref{alg:ILS} refer to an optional procedure, called online fine-tuning, to address this issue.
If the procedure is active, it first ensures that the estimated routing cost $V^*(s_0)$ of the selected solution accurately reflects the performance of the routing policy. This involves simulating the routing policy \eqref{eq:qpolicy} for a number of demand scenarios and confronting the observed routing costs with $V^*(s_0)$. We then perform a fine-tuning step, which consists in a partial re-training of the deep Q-Network. The main purpose is to locally improve the estimate of the value function. 

The procedure $\texttt{fine\_tuning}(x,f^r)$ takes as input the outsourcing solution $x$ and its simulated routing costs $f^r$. It then creates a copy of the pre-trained Q-network, denoted by $Q_{\theta'}(s_k, y)$, and updates its parameters $\theta'$ using the simulated costs, while keeping the original network $Q_{\theta}(s_k, y)$ fixed for routing committed customers.
The fine-tuning step is developed as a supervised learning procedure. 
Each call to $\texttt{fine\_tuning}(x, f^r)$ adds the pair $(x, f^r)$ to a training set. The model is then trained to minimize the error between the predicted routing cost $Q_{\theta'}(s_0, y_0)$ and the corresponding simulated cost $f^r$, where $y_0 = \arg\min_y Q_{\theta'}(s_0, y)$. The parameters are updated via backpropagation for a fixed number of epochs (10 in our implementation).

It is important to note that repeatedly simulating policy \eqref{eq:qpolicy} is time-consuming and cannot be executed every time the ILS requires evaluating routing costs. In consequence, the function  $\texttt{best\_improvement()}$, where the ILS algorithm explores the neighborhood of the current solution, completely relies on the cost estimation provided by the value function $V^*(s_k)$.
We use the simulation-based estimation only at the end of every round of neighborhood search.
We demonstrate the impact of adopting the fine-tuning step in Section \ref{sec:resultsidqncot}.

\section{Computational Results}\label{sec:experiments}
In this section we assess the performance of the proposed solution method through a set of computational experiments. This section is organized as follows. Section \ref{sec:expsetup} describes the experimental setup, namely the instances and the test and training protocol. Section \ref{sec:res1} reports the computational experiments, in which we assess the ILS algorithm and the routing policy \dqnco{} both in isolation and as components of the complete solution method, comparing each against custom benchmarks. Finally, Section \ref{sec:managerial} draws managerial insights from the results.

\subsection{Experimental setup}\label{sec:expsetup}
Section \ref{sec:instances} describes the instances used in our experiments and the procedure by which they are generated. Section \ref{sec:protocol} then specifies the test and training protocol, including the hyperparameters of the learning-based methods.

\subsubsection{Instances} \label{sec:instances}
Since no instance set for the \vrp{} is available in the literature, we constructed one. To this end, we built on the instance generation procedure of \cite{Dastpak2023}, proposed for a related problem, and adapted it to our setting. The problem studied in \cite{Dastpak2023}, denoted \vrpvcsd{}, has a hard duration limit constraint, so that no overtime operation is allowed, its objective is to maximize the demand served within that limit, and unserved customers are treated as missed opportunities rather than outsourced. Despite these differences, \vrpvcsd{} and \vrpx{} share the following features: variable customer sets, stochastic demands, and a fixed capacitated fleet operating under a duration limit. These commonalities allow us to adopt their generation procedure and extend it to account for the overtime and outsourcing costs of the \vrp{}.

Following \cite{Dastpak2023}, each instance $i$ is defined by a customer density $D\in\mathcal{D}=\{$Low, Moderate, High$\}$ and a vehicle capacity $Q\in\mathcal{Q}=\{25, 50, 75\}$. A customer density $D$ corresponds to a distribution function $\Gamma^C$ that governs the number of customers in the service area, their locations, drawn uniformly over that area, and their expected demands, each assigned uniformly at random from $\{5, 10, 15\}$. The average number of customers $\bar{n}$, the number of vehicles $m$, and the duration limit $L$ associated with each density $D$ are reported in Table \ref{tab:vrpvcsdinstances}; we refer the reader to \cite{Dastpak2023} for a complete description of the procedure.
In Section~\ref{sec:GAN}, we defined and used the parameter $n_{\max}$ as an upper bound for the number of customers in a given day. In our experiment we set $n_{\max} = 1.25 \bar{n}$.
The demand of each customer $c$ follows a uniform distribution $\Gamma^D_c$ on $[\bar{d}_c - \vartheta, \bar{d}_c + \vartheta]$, with $\vartheta=5$ for $\bar{d}_c\in\{10,15\}$ and $\vartheta=4$ otherwise.

\begin{table}[!htbp]
    \centering
    \footnotesize
    \begin{tabular}{c|c|c|c}
        $D$ & $\bar{n}$ & $m$ & $L$ \\ \hline        
        Low (L) & 23 & 3 & 221.47 \\
        Moderate (M) & 53 & 7 & 195.54 \\
        High (H) & 83 & 11 & 187.29 \\ \hline
    \end{tabular}
    \caption{Values of $\bar{n}, m,$ and $L$ for each density level $D$}
    \label{tab:vrpvcsdinstances}
\end{table}

\begin{figure}[!ht]
    \centering
    \includegraphics[width=0.3\textwidth]{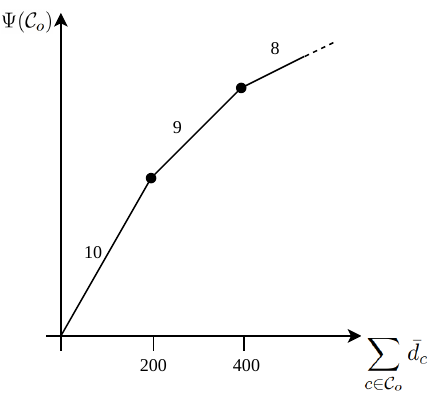}
    \qquad
    \begin{tabular}[b]{lr}\hline
      Quantity interval & Unit rate \\ \hline
      0-200 & 10.0 \\
      200-400 & 9.0 \\
      $400\leq$ & 8.0 \\ \hline
      & \\
      & \\
    \end{tabular}
    \caption{Outsourcing cost function per unit of expected total outsourced demand}
    \label{fig:outsourcingcost}
\end{figure}

Following \cite{Mendoza2016}, we set the overtime penalty to $\phi=2$. For the outsourcing cost function $\Psi(.)$, we adopt, as in \cite{Dabia2019}, a piecewise linear function of cost per unit of  expected total outsourced demand, illustrated in Figure \ref{fig:outsourcingcost}. Specifically, the numbers on each piece of this function represent the outsourcing cost rate. For example, if the planner decides to outsource a total of 300 demand volume, the outsourcing cost will be $200*10.0 + (300-200)*9.0=2900$. Under this function the LSP receives a discount when it outsources a larger expected volume of customer demand. The specific outsourcing cost rates are calibrated through a series of preliminary experiments on our instances, so as to balance outsourcing costs against routing costs, the latter comprising both travel and overtime costs. The rates must not be too low, as outsourcing costs well below routing costs would make outsourcing every customer the dominant decision for the LSP; conversely, they must not be too high, so that outsourcing part of the demand remains a viable option. In Section \ref{sec:managerial}, we further examine scenarios with lower and higher rates to analyze their impact on the performance of the solution.

In total, this yields nine instances, one for each pair $i\in\mathcal{D}\times\mathcal{Q}$, on which we compared our solution method against the benchmarks. 

\subsubsection{Test and training protocol}\label{sec:protocol}

We evaluate all methods on a common test set. For each instance $i$, a realization $\hat{i}$ is obtained by sampling a customer set from $\Gamma^C$ and, conditional on it, the customer demands from $\Gamma^D$. Each method was tested on every instance $i$ over $250{,}000$ realizations, formed by $500$ customer sets sampled from $\Gamma^C$ and, for each such set, $500$ demand realizations sampled from $\Gamma^D$.
In the ILS algorithm, the maximum number of iterations $\xi$ is set to $100$.
In the fine-tuning step, we use 50 randomly generated demand scenarios to re-evaluate the performance of the routing policy via simulation.

The Q-learning based methods require a training phase to develop their routing policy. 
Normally, the training set would be drawn from the marginal distribution $\Gamma^{\bar{\mathcal{C}}^o}$ of the committed set $\bar{\mathcal{C}}^o$. This is however not feasible in practice, since $\Gamma^{\bar{\mathcal{C}}^o}$ is induced by the outsourcing decisions searched by the ILS algorithm (as discussed in Section~\ref{sec:generalizedMDP}). 
We consequently construct the training set by adopting a surrogate procedure. Specifically, a realization $\hat{i}$ in the training set is obtained as follows. We first sample a customer set $\mathcal{C}$ from $\Gamma^C$. We then emulate the outsourcing decisions searched by the ILS algorithm by removing $\hat{n}$ customers at random from $\mathcal{C}$. The remaining customers are interpreted as the committed set $\bar{\mathcal{C}}^o$, thus forming the realization $\hat{i}$ . The value $\hat{n}$ is drawn uniformly from $[0, n-1]$. This procedure is repeated for all realizations in the training set.
We note that test instances and a code to generate train instances are available at \url{https://github.com/moda707/vrp-sdo}.

We now specify the hyperparameters of \dqnco{}. In Q-learning, the discount factor $\gamma$ is commonly set slightly below one to speed up convergence; however, as discussed earlier, $\gamma<1$ discounts the value function and thereby underestimates the routing costs. We therefore set $\gamma=1.0$. To further aid convergence, we employ a double Q-network, in which the target network is replaced by a copy of the primary network every $P^d=1000$ trials, together with a replay memory of size $50{,}000$. At each decision epoch, the network parameters $\theta$ are updated with probability $P^t=0.05$, and the minibatch size $|\tilde{B}|$ and embedding size $e$ are set to $32$ and $128$, respectively.

We trained \dqnco{} on five million realizations $\hat{i}$; each such realization used in training is referred to as a trial. The exploration rate $\epsilon$ decays linearly from $1.0$ to $0.1$ over the first third of the training trials, from $0.1$ to $0.05$ over the second third, and remains at $0.05$ thereafter, while the learning rate decays linearly from $10^{-3}$ to $10^{-4}$ over the first third.

All Q-learning-based methods are trained on GPUs, whereas all other components of method are executed on CPUs.
Specifically, the former were trained on Compute Canada GPU clusters equipped with NVIDIA Tesla V100 GPUs, while the latter were performed on the Béluga cluster of Compute Canada using one Intel Xeon Gold 6148 CPU core and 32 GB of RAM. 
Computational costs and training times are discussed in Section~\ref{sec:managerial}.

\subsection{Performance Analysis}\label{sec:res1}
To the best of our knowledge, no solution method in the literature can be applied directly to the \vrp{}. 
We therefore assess our method against a set of suitable custom benchmarks.
More specifically, we analyze the performance of the proposed solution method in three steps, progressively combining its two levels. In Section \ref{sec:resultsils}, we isolate the first level by assessing the ILS algorithm against a full-enumeration approach. In Section \ref{sec:resultsdqnco}, we isolate the second level by comparing the routing policy \dqnco{} against a set of benchmark routing policies on the routing problem alone. Finally, in Section \ref{sec:resultsidqncot}, we evaluate the complete method against benchmarks obtained by combining the ILS algorithm with each of these routing policies, all solving the \vrp{}.

\subsubsection{Performance of the ILS algorithm}\label{sec:resultsils}

We first assess the efficiency of the ILS algorithm employed in \idqnco{}. To this end, we compare it against a benchmark that fully enumerates the possible outsourcing solutions and, as in our method, uses \dqnco{} to evaluate their routing costs, so that the two differ only in how the outsourcing problem is searched. Since this benchmark has $O(2^n)$ time complexity, its computation time grows exponentially with the number of customers; we therefore restrict this experiment to low-density instances.

We ran the full-enumeration benchmark on instances with $D=$ Low and $Q\in\{25, 50, 75\}$, each over $10$ sampled customer sets, imposing a runtime limit of $48$ hours per customer set.
Within this limit, the enumeration returned the optimal outsourcing solution in $29$ of the $30$ customer sets. In every case, the optimal solution coincided with the solution returned by \idqnco{}, at a considerably lower computational cost: the full-enumeration benchmark required more than four hours on average, whereas the ILS algorithm obtained the same solution in eight seconds on average. This experiment confirms the effectiveness of the ILS algorithm. 
Table \ref{tab:resenum} reports the results, where ``Enum." relates to the full enumeration results.

\begin{table}[!ht]
\footnotesize
\centering
\begin{tabular}{c|cc|c|c}
\multirow{2}{*}{$Q$} & \multicolumn{2}{c|}{Computation Time} & \multirow{2}{*}{Enum. found Opt. Sol.} & \multirow{2}{*}{Enum. \& ILS Sol. matched} \\
\cline{2-3}
 & \idqnco{} (sec) & Enum. (sec) & & \\
\hline
25 & 12 & 11820 & 9/10 & 9/9 \\
50 & 6  & 15407 & 10/10 & 10/10 \\
75 & 5  & 19927 & 10/10 & 10/10 \\
\hline
\textbf{Avg:} & \textbf{8} & \textbf{15718} & \textbf{29/30} & \textbf{29/29} \\
\hline
\end{tabular}
\caption{Comparing the performance of the ILS with the full enumeration benchmark on instances with Low density}
\label{tab:resenum}
\end{table}

\subsubsection{Performance of \dqnco{} in \vrpx{}}\label{sec:resultsdqnco}

We now isolate the second level and assess the routing policy \dqnco{} on the \vrpx{}, comparing it against four benchmark routing policies: random policy (RP), greedy policy (GP), hyper-greedy policy (HP), and \qnco{}.
In RP, the next location for the active vehicle in decision epoch $k$ is chosen uniformly at random from the set of available customers $\overline{\mathcal{C}}^o_k = \{c\in\overline{\mathcal{C}}^o|h_c=1\}$. In GP, the active vehicle is sent to the nearest available customer, that is, $y^{\bar{v}}_k= \texttt{y}^D_c$ with $c = \arg\min_{c\in\overline{\mathcal{C}}^o_k} \tau_{\bar{v}, c}$. In HP, the active vehicle is sent to the customer with the highest ratio of served demand to travel distance, that is, $y^{\bar{v}}_k= \texttt{y}^D_c$ with $c=\arg\max_{c\in\overline{\mathcal{C}}^o_k} \frac{\min(\hat{d}_c, q_{\bar{v}})}{\tau_{\bar{v}, c}}$, so that customers with larger demands and shorter distances are prioritized.
Under all three policies, the active vehicle returns to the depot whenever $|\overline{\mathcal{C}}^o_k|=0$ or $q_{\bar{v}} = 0$; staying at the depot, that is, early trip termination, and preventive restocking are not available.
The fourth benchmark routing policy is \qnco{}, which is the Q-learning method of \cite{Dastpak2023} adapted to our problem. Specifically, it differs from \dqnco{} in how the observation function is constructed: as detailed in Section \ref{sec:rl}, \dqnco{} relies on a more sophisticated (i.e., deeper) neural network architecture that can capture more relevant information about the active vehicle. Accordingly, \qnco{} uses a shallower network and does not employ the graph-attention observation function of \dqnco{}. \qnco{} shares the same exploration and learning rate schedules as \dqnco{} (Section \ref{sec:protocol}) but is trained on three million realizations instead of five million.

We solved the \vrpx{} for each instance $i$ under all five routing policies, assuming that no customer is outsourced. For each instance, we evaluated the policies on the $250{,}000$ realizations $\hat{i}$ generated as described in Section \ref{sec:protocol}.
Table \ref{tab:resrouting} reports the results. The first two columns identify the instance; the next five report the routing cost — travel plus overtime — under RP, GP, HP, \qnco{}, and \dqnco{}; and the final four columns report the performance gap between \dqnco{} and  benchmark policy $X$, computed as 
$$\%X = \frac{\dqnco{} - X}{X} \times 100.$$

The results show that \dqnco{} outperforms all four benchmark routing policies on every instance. On average, it achieves routing costs that are $66.0\%$, $29.6\%$, and $34.6\%$ lower than those of RP, GP, and HP, and $19.6\%$ lower than those of \qnco{}. This gap demonstrates the effectiveness of the GAT-based observation function adopted in \dqnco{}.

\begin{table}[!ht]
\footnotesize
\centering
\begin{tabular}{cc|ccccc|cccc}
$D$ & $Q$  & RP      & GP     & HP     & \qnco{}   & \dqnco{}    & \%RP     & \%GP     & \%HP     & \%\qnco{}   \\ \hline
\multirow{3}{*}{L} & 25 & 2859.3  & 1835.8 & 1937.7 & 1702.3 & 1343.0 & -53.0\% & -26.8\% & -30.7\% & -21.1\% \\
                   & 50 & 2183.8  & 1119.5 & 1242.5 & 1031.6 & 748.1  & -65.7\% & -33.2\% & -39.8\% & -27.5\% \\
                   & 75 & 1943.1  & 857.9  & 964.1  & 874.4  & 577.8  & -70.3\% & -32.6\% & -40.1\% & -33.9\% \\ \hline
\multirow{3}{*}{M} & 25 & 6987.6  & 4065.9 & 4264.9 & 3572.3 & 2945.1 & -57.9\% & -27.6\% & -30.9\% & -17.6\% \\
                   & 50 & 5394.5  & 2455.5 & 2703.3 & 2174.0 & 1660.6 & -69.2\% & -32.4\% & -38.6\% & -23.6\% \\
                   & 75 & 4831.9  & 1853.1 & 1987.9 & 1704.1 & 1343.7 & -72.2\% & -27.5\% & -32.4\% & -21.1\% \\ \hline
\multirow{3}{*}{H} & 25 & 11226.6 & 6122.1 & 6448.1 & 4967.4 & 4483.5 & -60.1\% & -26.8\% & -30.5\% & -9.7\%  \\
                   & 50 & 8707.8  & 3701.3 & 4025.0 & 2873.2 & 2298.5 & -73.6\% & -37.9\% & -42.9\% & -20.0\% \\
                   & 75 & 7832.6  & 2779.2 & 2928.8 & 2221.9 & 2181.0 & -72.2\% & -21.5\% & -25.5\% & -1.8\%  \\ \hline
\multicolumn{7}{r|}{\textbf{Avg:}} & \textbf{-66.0\%} & \textbf{-29.6\%} & \textbf{-34.6\%} & \textbf{-19.6\%} \\ \hline
                          
\end{tabular}
\caption{Comparing the performance of the routing policy \dqnco{} with the other four routing policies RP, GP, HP, and \qnco{}}
\label{tab:resrouting}
\end{table}

\subsubsection{Performance of the full solution method}\label{sec:resultsidqncot}
Having assessed the ILS algorithm and the routing policy \dqnco{} in isolation, we now evaluate the full solution method, in which the two levels operate jointly to solve the \vrp{}. We combine the ILS algorithm with each of the routing policies introduced in Section \ref{sec:resultsdqnco}, yielding the methods \irp{}, \igp{}, \ihp{}, \iqnco{}, and \idqnco{}. Applying the online fine-tuning step of Section \ref{sec:ilsenhancement} to the two Q-learning-based methods, \iqnco{} and \idqnco{}, yields \iqncot{} and \idqncot{}. 

We evaluated all seven methods on the test set described in Section \ref{sec:protocol}. For each of the $500$ customer sets, which defines a daily \vrp{}, every method produces an outsourcing solution by searching with the ILS algorithm. During this search, the ILS estimates the routing cost of each candidate solution as follows: for \irp{}, \igp{}, and \ihp{}, by simulating the corresponding policy over $50$ demand realizations on the committed customer set $\overline{\mathcal{C}}^o$ and averaging; for the Q-learning-based methods, by using the trained value function evaluated at the initial post-decision state, $V^*(s_0)$. We then evaluate the chosen outsourcing solution by simulating the method's routing policy on the committed customer set over $500$ demand realizations, computing the outsourcing cost via the function $\Psi(\cdot)$ introduced in Section \ref{sec:instances}, and reporting their sum as the total cost.

Table \ref{tab:res1} reports the total cost of each method on every instance, where the total cost is the sum of the outsourcing cost and the estimated routing cost. We also report performance gaps relative to \idqncot{}. For each method $X\in\{\text{\irp{}}, \text{\igp{}}, \text{\ihp{}}, \text{\iqnco{}}, \text{\iqncot{}}, \text{\idqnco{}}\}$, the gap is computed as
\[
\%X = \frac{\text{\idqncot{}} - X}{X} \times 100.
\]

The results show that \idqncot{} achieves the lowest total cost on every instance, with average gaps of $-45.3\%$, $-22.3\%$, $-27.1\%$, $-19.0\%$, $-13.7\%$, and $-1.2\%$ over \irp{}, \igp{}, \ihp{}, \iqnco{}, \iqncot{}, and \idqnco{}, respectively. The gap between \idqnco{} and \iqnco{}, the two methods that differ only in the routing policy used, confirms the effect of the GAT-based observation function already observed in Section \ref{sec:resultsdqnco} for \dqnco{} on \vrpx{}. Online fine-tuning improves both \iqnco{} and \idqnco{}, more substantially for the shallower one: an average gap of $-6.3\%$
for \iqncot{} over \iqnco{}, compared to $-1.2\%$ for \idqncot{} over \idqnco{}.
We examine the fine-tuning step effects on the first-level decisions separately in the next experiment. A potential interpretation is that the more accurate value-function estimates of \dqnco{} leave less residual error for the fine-tuning step to correct.

\begin{table}[!ht]
\footnotesize
\centering
\setlength{\tabcolsep}{0.75ex}
\begin{tabular}{cc|ccccccc|cccccc}
$D$ & $Q$  & \irp{}     & \igp{}      & \ihp{}      & \iqnco{} & \iqncot{}    & \idqnco{} & \textbf{\idqncot{}}   & \%\irp{}     & \%\igp{}     & \%\ihp{}     & \%\iqnco{} & \%\iqncot{} & \%\idqnco{}   \\ \hline
\multirow{3}{*}{L} & 25 & 1675.3 & 1520.8 & 1548.0 & 1475.9 & 1371.1 & 1325.6 & 1278.6  & -23.7\% & -15.9\% & -17.4\% & -13.4\% & -6.7\%   & -3.5\%   \\
                   & 50 & 1429.9 & 1003.3 & 1096.7 & 1000.2 & 878.2  & 742.3  & 737.7   & -48.4\% & -26.5\% & -32.7\% & -26.2\% & -16.0\%  & -0.6\%   \\
                   & 75 & 1289.0 & 753.3  & 804.1  & 867.2  & 732.1  & 577.1  & 572.3   & -55.6\% & -24.0\% & -28.8\% & -34.0\% & -21.8\%  & -0.8\%   \\ \hline
\multirow{3}{*}{M} & 25 & 3962.5 & 3584.9 & 3629.5 & 3406.5 & 3322.5 & 2939.6 & 2913.1  & -26.5\% & -18.7\% & -19.7\% & -14.5\% & -12.3\%  & -0.9\%   \\
                   & 50 & 3435.3 & 2263.8 & 2549.7 & 2172.1 & 2057.6 & 1661.4 & 1643.8  & -52.1\% & -27.4\% & -35.5\% & -24.3\% & -20.1\%  & -1.1\%   \\
                   & 75 & 3163.5 & 1629.4 & 1721.4 & 1707.4 & 1587.2 & 1338.2 & 1316.4  & -58.4\% & -19.2\% & -23.5\% & -22.9\% & -17.1\%  & -1.6\%   \\ \hline
\multirow{3}{*}{H} & 25 & 6136.5 & 5617.9 & 5710.5 & 4971.4 & 4932.2 & 4482.8 & 4478.5  & -27.0\% & -20.3\% & -21.6\% & -9.9\%  & -9.2\%   & -0.1\%   \\
                   & 50 & 5399.3 & 3589.7 & 4394.6 & 2871.0 & 2825.3 & 2304.1 & 2298.6  & -57.4\% & -36.0\% & -47.7\% & -19.9\% & -18.6\%  & -0.2\%   \\
                   & 75 & 5050.0 & 2409.4 & 2524.0 & 2227.5 & 2123.2 & 2130.0 & 2095.8  & -58.5\% & -13.0\% & -17.0\% & -5.9\%  & -1.3\%   & -1.6\%   \\ \hline
\multicolumn{9}{r|}{\textbf{Avg:}}   & \textbf{-45.3\% }& \textbf{-22.3\%} & \textbf{-27.1\%} & \textbf{-19.0\%} & \textbf{-13.7\%}  & \textbf{-1.2\%}   \\ \hline
\end{tabular}
\caption{Comparing the performance of the \idqncot{} with the other six methods in terms of the total costs}
\label{tab:res1}
\end{table}

\subsection{Managerial insights}\label{sec:managerial}
The experiments of Section \ref{sec:res1} compare the solution methods against one another. We now examine two questions that an LSP would face when deploying these methods in practice. First, how sensitive are the resulting decisions to the outsourcing rates charged by the common carrier, and how does this sensitivity depend on the LSP's fleet capacity? Second, how long does each method take to run, and how does that constrain the choice of method given the LSP's daily decision window?

We consider four scenarios for the outsourcing rates charged by the common carrier, denoted Very Low, Low, Moderate, and High, with unit rates equal to $0.2$, $0.5$, $1$, and $2$ times the rates used in Section \ref{sec:res1}, respectively. For each scenario, we ran \idqncot{} on the three instances with Moderate density and $Q\in\{25, 50, 75\}$, on the first $50$ customer sets of the test set described in Section \ref{sec:protocol}. Table \ref{tab:outsourcinganalysis} reports, for each scenario--instance pair, the percentage of outsourced customers (\% Out.\ Cust.), the percentage of travel time spent in overtime (\%Overtime), the routing cost (RC), the outsourcing cost (OC), the total cost (TC), and the percentage of total cost paid for outsourcing (\% OC/TC).

Across all capacities, higher rates lead to higher total costs, fewer outsourced customers, and a lower percentage of total cost spent on outsourcing. The intensity of this response, however, depends on the vehicle capacity. As rates rise from Very Low to High, the percentage of outsourced customers drops from $96.5\%$ to $3.0\%$ at $Q=25$, but only from $26.5\%$ to $1.1\%$ at $Q=75$. An LSP with a fleet composed of  vehicles with a small capacity is therefore exposed to the carrier's pricing: most of its operation runs through the carrier at low rates and almost none at high rates, so its operational mode heavily depends on the carrier's price level. An LSP with a fleet composed of large-capacity vehicles, in contrast, is largely rate-insensitive, since even at the lowest rates only about a quarter of demand is outsourced and at the highest rates outsourcing is essentially abandoned.

We now turn to the computation time required to obtain the outsourcing decision each day, and how it constrains the choice of method.
We measured the average CPU time required by each method to produce the outsourcing solution on each instance, with an imposed time limit of two hours per run. This time limit was chosen such that, on average, more than 90\% of the experiments (57 out of 63, corresponding to 9 instances and 7 methods) were able to find an outsourcing solution within the given time limit. Table \ref{tab:rescomp} reports the results. The three simulation-based methods, \irp{}, \igp{}, and \ihp{}, are considerably slower than the learned ones: on the high-density instances the time limit is reached on almost every run, and on average across all instances they require more than an hour each. \iqnco{} and \idqnco{}, in contrast, complete in $18$ and $43$ seconds on average, respectively, with the difference reflecting the deeper network used by \dqnco{} and therefore the higher per-call cost of estimating routing costs. Applying the online fine-tuning step raises these times to $3.2$ minutes for \iqncot{} and $8.5$ minutes for \idqncot{}. These daily times do not include the offline training required by the Q-learning-based methods. Training \qnco{} takes between one and two days per instance $i\in\mathcal{D}\times\mathcal{Q}$, and training \dqnco{} between four and sixteen days; both were trained on Compute Canada GPUs. This time is spent once before deployment and does not recur.

Generally speaking, simulation-based methods, i.e., \irp{}, \igp{}, and \ihp{}, rely on online simulation and become slower as the problem size grows, which may render them unsuitable when the LSP's daily decision window is limited. Learning-based methods avoid this online cost by transferring computation effort to the offline training phase, and consequently appear better suited to the operational needs of an LSP, at the price of a substantial one-time training investment that they require beforehand.

\begin{table}[!ht]
    \scriptsize
    \centering
    \setlength{\tabcolsep}{0.75ex}
    \begin{tabular}{c|c|cccccc}
    $\Psi$ & $Q$  & \% Out. Cust. & \%Overtime & RC & OC & TC & \% OC/TC \\ \hline
    \multirow{3}{*}{{[}2, 1.8, 1.6{]}}  & 25 & 96.51\%                 & 1.09\%     & 44.9          & 922.0             & 966.9       & 95.36\%              \\
 & 50 & 53.36\%                 & 7.70\%     & 659.8         & 498.2             & 1158.0      & 43.02\%              \\
 & 75 & 26.45\%                 & 5.78\%     & 861.0         & 243.4             & 1104.4      & 22.04\%              \\ \hline
    \multicolumn{2}{r|}{\textbf{Avg:}}                 & 58.77\%                 & 6.44\%     & 521.9         & 554.5             & 1076.5      & 51.52\%              \\ \hline
    \multirow{3}{*}{{[}5, 4.5, 4{]}}    & 25 & 47.02\%                 & 15.72\%    & 1365.2        & 1004.7            & 2369.9      & 42.39\%              \\
 & 50 & 8.68\%                  & 17.58\%    & 1431.2        & 150.5             & 1581.7      & 9.51\%               \\
 & 75 & 6.07\%                  & 8.77\%     & 1150.8        & 116.0             & 1266.8      & 9.16\%               \\ \hline
    \multicolumn{2}{r|}{\textbf{Avg:}}                 & 20.59\%                 & 14.25\%    & 1315.7        & 423.7             & 1739.5      & 24.36\%              \\ \hline
    \multirow{3}{*}{\textbf{{[}10, 9, 8{]}}}    & 25 & 11.06\%                 & 31.89\%    & 2517.2        & 377.0             & 2894.2      & 13.03\%              \\
 & 50 & 3.26\%                  & 19.03\%    & 1540.4        & 97.0              & 1637.4      & 5.92\%               \\
 & 75 & 2.23\%                  & 11.54\%    & 1252.3        & 66.0              & 1318.3      & 5.01\%               \\ \hline
    \multicolumn{2}{r|}{\textbf{Avg:}}                 & 5.52\%                  & 22.76\%    & 1770.0        & 180.0             & 1950.0      & 9.23\%               \\ \hline
    \multirow{3}{*}{{[}15, 13.5, 12{]}} & 25 & 2.99\%                  & 35.80\%    & 2801.8        & 144.0             & 2945.8      & 4.89\%               \\
 & 50 & 1.54\%                  & 20.09\%    & 1583.2        & 61.5              & 1644.7      & 3.74\%               \\
 & 75 & 1.11\%                  & 11.31\%    & 1272.8        & 45.0              & 1317.8      & 3.41\%               \\ \hline
    \multicolumn{2}{r|}{\textbf{Avg:}}                 & 1.88\%                  & 25.03\%    & 1885.9        & 83.5              & 1969.4      & 4.24\%       \\ \hline       
    \end{tabular}
    \caption{Sensitivity analysis on the outsourcing costs structure}
\label{tab:outsourcinganalysis}
\end{table}

\begin{table}[!ht]
\footnotesize
\centering
\setlength{\tabcolsep}{0.75ex}
\begin{tabular}{cc|ccccccc}
$D$ & $Q$  & \irp{}   & \igp{}   & \ihp{}   & \iqnco{} & \iqncot{} & \idqnco{} & \idqncot{} \\ \hline
\multirow{3}{*}{L}      & 25 & 276  & 328  & 473  & 5    & 52    & 12    & 140    \\
                          & 50 & 265  & 236  & 382  & 4    & 34    & 6     & 78     \\
                          & 75 & 225  & 303  & 299  & 3    & 38    & 5     & 73     \\ \hline
\multicolumn{2}{r|}{\textbf{Avg:}} & 255  & 289  & 385  & 4    & 41    & 8     & 97     \\ \hline
\multirow{3}{*}{M} & 25 & 3566 & 6215 & 6698 & 31   & 204   & 95    & 644    \\
                          & 50 & 4167 & 2858 & 5191 & 9    & 109   & 30    & 343    \\
                          & 75 & 3442 & 3101 & 3831 & 7    & 183   & 18    & 424    \\ \hline
\multicolumn{2}{r|}{\textbf{Avg:}}  & 3725 & 4058 & 5240 & 16   & 166   & 48    & 471    \\ \hline
\multirow{3}{*}{H}     & 25 & 7200 & 7200 & 7200 & 72   & 684   & 124   & 1055   \\
                          & 50 & 7200 & 6919 & 7200 & 17   & 221   & 44    & 923    \\
                          & 75 & 7200 & 7033 & 7095 & 13   & 192   & 56    & 905    \\ \hline
\multicolumn{2}{r|}{\textbf{Avg:}}       & 7200 & 7051 & 7165 & 34   & 366   & 75    & 961    \\ \hline
\multicolumn{2}{r|}{\textbf{Total Avg:}} & 3727 & 3799 & 4263 & 18   & 191   & 43    & 510       \\ \hline
\end{tabular}
\caption{Comparing the computation time (in Seconds) required to solve each instance}
\label{tab:rescomp}
\end{table}

\section{Conclusions}\label{sec:conclusion}
In this study, we addressed a stochastic variant of the VRP with outsourcing options, in which a set of customers along with their demand distributions are revealed at the beginning of each day. Given this set, an LSP must determine the committed and outsourcing sets of customers. The LSP may pay overtime costs in case of violating a daily work shift limit. The objective is to determine the set of outsourced customers and find the optimal routing policy to serve committed customers such that travel, overtime, and outsourcing costs are minimized.

We formalized this problem as a two-level problem, where the first level identifies committed and outsourced customers, while the second level estimates routing costs. 
To solve this problem, we introduce the \idqnco{}, a heuristic algorithm featuring an ILS algorithm. It evaluates various customer partitions (i.e., committed vs outsourced) and solves a dynamic VRP with stochastic demands for each explored partition to compute its routing costs. To efficiently handle the routing costs estimation, we formalized the second level routing problem as an MDP and proposed  the \dqnco{}, which is a DQN-based algorithm, to solve it. 
The \dqnco{} required an intensive training phase, which was in conflict with the need to solve the \vrp{} on a daily basis. To overcome this limitation, we modified the MDP formulation of the routing problem in such a way as to enable the \dqnco{} to handle random customer sets drawn from a suitable probability distribution. By doing so, our \dqnco{} was trained offline once, allowing the \idqnco{} to be readily used on a daily basis. 
Additionally, we proposed \idqncot{}, which enhanced \idqnco{} with an online fine-tuning step that improved the accuracy of routing costs estimation as the ILS algorithm explores different customer partitions.

We compared \idqncot{} against six benchmarks, all employing the same ILS algorithm but differing in their routing policy and, for Q-learning-based methods, in whether the online fine-tuning step was applied. 
The results revealed that \idqncot{} consistently outperformed \idqnco{} and all other benchmarks, reducing overall costs by more than 22.3\% relative to benchmarks with basic routing policies and by 13.7\% on average relative to \iqncot{}. 
Notably, \idqncot{} demonstrated superior performance compared to \iqnco{}, illustrating the effectiveness of the proposed GAT structure introduced in \dqnco{}. Additionally, the superiority of \idqncot{} over \idqnco{} shows the effectiveness of online fine-tuning in improving the accuracy of the pre-trained value function for estimating routing costs.

From a managerial standpoint, our results offer practical guidance for LSPs considering the \vrp{}  operational setting. The sensitivity analysis on outsourcing rates shows that the impact of a common carrier's pricing on an LSP's operations depends heavily on its fleet composition.  LSPs operating small capacity vehicles see their outsourced share shift dramatically with rate changes, while LSPs operating large capacity vehicles remain comparatively insensitive to the carrier's pricing.
This suggests that fleet sizing decisions should account not only for expected demand, but also for the LSP's exposure to third-party pricing dynamics.
Furthermore, our results show that the DQN-based methods generate outsourcing decisions in under a minute on average, well within an LSP's daily decision window, whereas classical approaches are often too slow for this purpose, sometimes taking up to two hours.
This computational advantage comes at the cost of an offline training phase of several days, which is incurred once before deployment and does not recur, making the proposed methodology particularly attractive for LSPs seeking a scalable, reusable solution to a problem they must solve every day.

This study demonstrates that an offline-trained Q-network can successfully provide nearly instantaneous and accurate estimations of second-level problems. These estimations can be integrated into traditional heuristic schemes, such as an ILS in our case, which explore the first-level decision space. Therefore, it is promising to explore the application of this framework to other complex optimization problems, particularly those where the second-level problem is a stochastic dynamic problem.

\bibliography{references}
\bibliographystyle{elsarticle-harv}

\appendix
\section{Table of notations}
\begin{longtable}{l|l}
\textbf{Symbol} & \textbf{Description} \\
\hline
\endfirsthead

\multicolumn{2}{l}{\textit{Table of Notations (continued)}}\\
\textbf{Symbol} & \textbf{Description} \\
\hline
\endhead

\hline
\multicolumn{2}{r}{\textit{continued on next page}}\\
\endfoot

\hline
\caption{Table of Notations}\\
\endlastfoot

\multicolumn{2}{l}{\textbf{Problem}} \\ \hline
$c, \mathcal{C}, \overline{\mathcal{C}}^o, \mathcal{C}^o$ & Customer $c$ and the set of realized, committed and outsourced customers \\
$v, \mathcal{V}$ & Vehicle $v$ and the set of vehicles \\
$G, N, E$ & A complete graph $G$ with nodes $N$ and arcs $E$ \\
$l_0, l_c, l_v$ & Location of the depot, customer $c$, and vehicle $v$ \\
$\tau_{ij}$ & Travel time between locations $i$ and $j$ \\
$h_c, g_v$ & Availability of customer $c$ and operational status of vehicle $v$\\
$\bar{d}_c, d_c, \hat{d}_c$ & Expected, realized, and unserved demand of customer $c$ \\
$\Gamma^D_c, ~w$ & Probability distribution of demand for customer $c$ and a vector of demand realizations \\
$Q, q_v$ & Total capacity of each vehicle and the remaining capacity of vehicle $v$ \\
$a_v, L$ & Arrival time of vehicle $v$ at its next location and the shift duration of each vehicle \\
$\Psi(\mathcal{C}^o), \phi$ & Cost function to outsource customers $\mathcal{C}^o$ and the overtime cost penalty \\
$\pi_r,\pi^*_r$ & A routing policy and the optimal routing policy \\
$T_{\pi_r}(\bar{\mathcal{C}}^o, v, w)$ & Duration of the route performed by $v$, given $\bar{\mathcal{C}}^o$, $w$ and $\pi_r$ \\
$R(\bar{\mathcal{C}}^o)$ & Minimum expected routing costs for serving $\bar{\mathcal{C}}^o$ \\
$x_c$ & A binary decision variable for the first level problem \\

\hline
\multicolumn{2}{l}{\textbf{MDP}} \\ \hline
$k, ~t_k, ~s_k$ & Decision epoch, its time, and the system state \\
$O(s_k, \bar{v}), o_{k,\bar{v}}$ & Observation function and observation of active vehicle $\bar{v}$ at decision epoch $k$\\
$F^{\mathcal{C}}, F^{\mathcal{N}}, F^{\mathcal{V}}$ & State of customers, all nodes (customers + depot), and vehicles \\
$\bar{\mathcal{V}}_k$, $\bar{v}$ & Set of vehicles active at decision epoch $k$ and the active vehicle \\
$\overline{\mathcal{C}}^o_k, ~ \tilde{\mathcal{C}}_k$ & Set of available customers and set of customers being served at decision epoch $k$ \\
$y_k$ & Action vector at decision epoch $k$ \\
$y^D_c$, $y^I_c$ & Direct and indirect visit actions to customer $c$ \\
$A(s_k),\bar{A}(s_k)$ & Action space at state $s_k$ for MDP and MDP-CO \\
$\gamma$ & Discount factor in the MDP \\
$C(s_k, y_k)$ & Cost of taking action $y_k$ in state $s_k$ \\
$\Gamma^C, \Gamma^{\overline{\mathcal{C}}^o}$ & Distribution of daily customer realizations and committed customer subsets $\overline{\mathcal{C}}^o$ \\
$S_{\overline{\mathcal{C}}^o}$, $\bar{S}$, $\Omega$ & State space for a fixed customer subset, generalized state space, and active vehicle's observation space \\
$V^{\pi_r}(s_k), V^{\pi_r}(o_{k,\bar{v}})$ & Value function under routing policy $\pi_r$ for both state and observation\\
$Q^{\pi_r}(s_k,y_k), Q^{\pi_r}(o_{k,\bar{v}}, y_k)$ & Q-factors under routing policy $\pi_r$ for both state and observation\\
$\hat{Q}(s_k, y_k)$ & Estimated expected Q-factor of state $s_k$ and action $y_k$\\
$V^*(s_k), V^*(o_{k,\bar{v}})$ & Optimal value function for a given state $s_k$ and observation $o_{k,\bar{v}}$ \\
$Q^*(s_k,y_k), Q^*(o_{k,\bar{v}}, y_k)$ & Q-factors of action $y_k$ for a given state $s_k$ and observation $o_{k,\bar{v}}$\\

\hline
\multicolumn{2}{l}{\textbf{Neural network and ILS}} \\ \hline
$e$ & Embedding dimension \\
$H^{\mathcal{N}}$, $H^{\mathcal{V}}$ & Nodes and vehicles embeddings \\
$H^{\bar{v}}, H^{\mathcal{G}}$ & Embeddings of the active vehicle and the graph \\
$Att(F, \bar{F})$ & Attention mechanism on query $F$ and context $\bar{F}$ \\
$Q_\theta(s, y)$ & Q-factor function for state $s$, action $y$, and deep neural network parameters of $\theta$ \\
$\theta, \bar{\theta}$ & Parameters of the primary and target deep neural networks \\
$B$ & Memory replay buffer \\
$P^d, P^t$ & Double Q-network replacement interval and Q-factors update probability at each decision epoch \\
$\xi$ & Maximum allowed iterations for the ILS algorithm \\

\hline
\multicolumn{2}{l}{\textbf{Instances}} \\ \hline
$\mathcal{D}, \mathcal{Q}$ & Sets of customer densities and vehicle capacities \\

\end{longtable}

\end{document}